\documentclass[11pt]{amsart}
\usepackage{amssymb}
\usepackage{amsmath}
\usepackage{mathrsfs}
\usepackage{amsthm,enumerate}
\usepackage{graphicx}
\usepackage{multirow}
\usepackage{url}
\usepackage{subfigure,verbatim}
\usepackage{hyperref}
\usepackage{cite}
\usepackage{color}

\newcommand{\calH}{\mathcal{H}}
\newcommand{\calR}{\mathcal{R}}
\newcommand{\calF}{\mathcal{F}}

\DeclareMathOperator{\argmin}{argmin}

\newcommand{\E}{\mathbb{E}}
\newcommand{\R}{\mathbb{R}}
\newcommand{\abs}[1]{\left |#1\right |}
\newcommand{\norm}[1]{\left \|#1\right\|}

\newcommand{\paren}[1]{\left(#1\right)}
\newcommand{\bracket}[1]{\left[#1\right]}
\newcommand{\set}[1]{\left\{#1\right\}}
\newcommand{\inner}[2]{\left\langle #1,#2\right\rangle}
\newcommand{\eps}{\epsilon}

\newtheorem{lem}{Lemma}[section]
\newtheorem{cor}{Corollary}[section]
\newtheorem{prop}{Proposition}[section]
\newtheorem{thm}{Theorem}[section]

\newtheorem{assume}{Assumption}[]

\numberwithin{equation}{section}

\title[Relative Entropy \& Robbins-Monro]{Relative Entropy
  Minimization over Hilbert Spaces via Robbins-Monro}
\author{G. Simpson}
\email{grs53@drexel.edu}
\address{Department of Mathematics, Drexel University, Philadelphia,
  PA 19104, USA}
\author{D. Watkins}
\address{Department of Ocean, Earth, and Atmospheric Sciences, Oregon State University, Corvallis, OR 97331, USA}


\begin{document}

\begin{abstract}
  One way of getting insight into non-Gaussian measures is to first obtain good
Gaussian approximations.  These best fit Gaussians can then provide a sense of
the mean and variance of the distribution of interest. They can also be used to
accelerate sampling algorithms.  This begs the question of how one should
measure optimality, and how to then obtain this optimal approximation.  Here, we
consider the problem of minimizing the distance between a family of Gaussians
and the target measure with respect to relative entropy, or Kullback-Leibler
divergence.  As we are interested in applications in the infinite dimensional
setting, it is desirable to have convergent algorithms that are well posed on
abstract Hilbert spaces. We examine this minimization problem by seeking roots
of the first variation of relative entropy, taken with respect to the mean of
the Gaussian, leaving the covariance fixed. We prove the convergence of
Robbins-Monro type root finding algorithms in this context, highlighting the
assumptions necessary for convergence to relative entropy minimizers.  Numerical
examples are included to illustrate the algorithms.
\end{abstract}

\maketitle

\section{Introduction}

In \cite{Pinski:2015jn,Pinski:2015jy,
Lu_Gaussian1_2017,Lu_Gaussian2_2017,sanz2017gaussian}, it was proposed that
insight into a probability distribution, $\mu$, posed on a Hilbert space,
$\calH$, could be obtained by finding a best fit Gaussian approximation, $\nu$.
This notion of best, or optimal, was with respect to the relative entropy,
or Kullback-Leibler divergence:
\begin{equation}
  \label{e:kl}
  \calR(\nu||\mu) =\begin{cases}
    \E^{\nu}\bracket{\log \frac{d\nu}{d\mu}}, & \nu \ll \mu,\\
    +\infty, &\text{otherwise}.
  \end{cases}
\end{equation}
Having a Gaussian approximation provides qualitative insight into $\mu$, as it
provides a concrete notion of the mean and variance of the distribution.
Additionally, this optimized distribution can be used in algorithms, such as
random walk Metropolis, as a preconditioned proposal distribution to improve
performance.  Such a strategy can benefit a number of applications, including
path space sampling for molecular dynamics and parameter estimation in
statistical inverse problems.

Observe that in the definition of $\calR$, \eqref{e:kl}, there is an
asymmetry in the arguments.  Were we to work with $\calR(\mu||\nu)$,
our optimal Gaussian would capture the first and second moments of
$\mu$, and in some applications this is desirable.  However, for a
multimodal problem (consider a distribution with two well separated
modes), this would be inadequate; our form attempts to match
individual modes of the distribution by a Gaussian.  For a recent review of the $\calR(\nu||\mu)$ problem, see \cite{Blei:2017798}, where it is remarked that this choice of arguments is likely to underestimate the dispersion of the distribution of interest, $\mu$.  The other
ordering of arguments has been explored, in the finite dimensional case, in
\cite{Andrieu:2006hg,Andrieu:2008kh,Haario:2001gu,Roberts:2007ey}.

To be of computational use, it is necessary to have an algorithm that
will converge to this optimal distribution.  In \cite{Pinski:2015jy},
this was accomplished by first expressing $\nu = N(m, C(p))$, where
$m$ is the mean and $p$ is a parameter inducing a well defined
covariance operator, and then solving the problem,
\begin{equation}
  (m,p) \in \argmin \calR(N(m, C(p))||\mu),
\end{equation}
over an admissible set.  The optimization step itself was done using the
Robbins-Monro algorithm (RM), \cite{Robbins:1950ua}, by seeking a root of the
first variation of the relative entropy. While the numerical results of
\cite{Pinski:2015jy} were satisfactory, being consistent with theoretical
expectations, no rigorous justification for the application of RM to the
examples was given.

In this work, we emphasize the study and application of RM to potentially
infinite dimensional problems.  Indeed, following the framework of
\cite{Pinski:2015jn,Pinski:2015jy}, we assume that $\mu$ is posed on the Borel
$\sigma$-algebra of a separable Hilbert space $(\calH, \inner{\bullet}{\bullet},
\norm{\bullet})$. For simplicity, we will leave the covariance operator $C$
fixed, and only optimize over the mean, $m$.  Even in this case, we are seeking
$m\in \calH$, a potentially infinite-dimensional space.

\subsection{Robbins-Monro}

Given the objective function $f:\calH \to \calH$, assume that it has a
root, $x_\star$.  In our application to relative entropy, $f$ will be
its first variation. Further, we assume that we
can only observe a noisy version of $f$,
$F:\calH \times \chi\to \calH$, such that for all $x \in \calH$,
\begin{equation}
  f(x) = \E[F(x,Z)] = \int_{\chi} F(x, z) \mu_Z(dz),
\end{equation}
where $\mu_Z$ is the distribution associated with the random variable (r.v.)
$Z$, taking values in the auxiliary space $\chi$.  The naive Robbins-Monro
algorithm is given by
\begin{equation}
  \label{e:rm_naive}
  X_{n+1} = X_n - a_{n+1} F(X_n, Z_{n+1}),
\end{equation}
where $Z_n \sim \mu_Z$, are independent and identically distributed
(i.i.d.), and $a_n>0$  is a carefully chosen sequence.
Subject to assumptions on $f$, $F$, and the distribution $\mu_Z$, it
is known that $X_n$ will converge to $x_\star$ almost surely (a.s.),
in finite dimensions, \cite{Robbins:1950ua,Blum:1954tf,Blum:1954vm}.
Often, one needs to assume that $f$ grows at most
linearly,
\begin{equation}
  \label{e:linear_bound}
  \norm{f(x)}\leq c_0 + c_1 \norm{x},
\end{equation}
in order to apply the results in the aforementioned papers.
The analysis in the finite dimensional case has been refined
tremendously over the years, including an analysis based on continuous
dynamical systems. We refer the reader to the books
\cite{Kushner:2003aa,Albert:1967aa,chen2002stochastic} and references therein.

\subsection{Trust Regions and Truncations}

As noted, much of the analysis requires the regression function
$f$ to have, at most, linear growth.  Alternatively, an {\it a priori}
assumption is sometimes made that the entire sequence generated by
\eqref{e:rm_naive} stays in a bounded set.  Both assumptions are
limiting, though, in practice, one may find that the algorithms
converge.

One way of
overcoming these assumptions, while still ensuring convergence, is to introduce trust regions that the sequence $\{X_n\}$
is permitted to explore, along with a ``truncation'' which enforces
the constraint.  Such truncations distort \eqref{e:rm_naive} into
\begin{equation}
  \label{e:rm_truncation1}
  X_{n+1} = X_n - a_{n+1} F(X_n, Z_{n+1}) + a_{n_+1} P_{n+1},
\end{equation}
where $P_{n+1}$ is the projection keeping the sequence $\{X_n\}$ within the
trust region.  Projection algorithms are also discussed in
\cite{Kushner:2003aa,Albert:1967aa,chen2002stochastic}.

We consider RM on a possibly infinite dimensional separable Hilbert space.  This
is of particular interest as, in the context of relative entropy optimization,
we may be seeking a distribution in a Sobolev space associated with a PDE model.
A general analysis of RM with truncations in Hilbert spaces can be found in
\cite{Yin:1990wv}.  The main purpose of this work is to adapt the analysis of
\cite{Lelong:2008ck} to the Hilbert space setting for two versions of the
truncated problem.  The motivation for this is that the analysis of
\cite{Lelong:2008ck} is quite straightforward, and it is instructive to see how
it can be easily adapted to the infinite dimensional setting.  The key
modification in the proof is that results for Banach space valued martingales
must be invoked.  We also adapt the results to a version of the algorithm where
there is prior knowledge on the location of the root.  With these results in
hand, we can then verify that the relative entropy minimization problem can be
solved using RM.

\subsubsection{Fixed Trust Regions}

In some problems, one may have {\it a priori} information on the root.
For instance, we may know that $x_\star \in U_1$, some open bounded set.  In this version of the
truncated algorithm, we have two open bounded sets,
$U_0\subsetneq U_1$, and $x_\star \in U_1$.  Let $\sigma_0 = 0$ and
$X_0 \in U_0$ be given, then \eqref{e:rm_truncation1} can be
formulated as
\begin{subequations}
  \label{e:fixed_trust}
  \begin{gather}
    \label{e:proposal1}
    \tilde{X}_{n+1} = X_n - a_{n+1} F(X_n, Z_{n+1})\\
    \label{e:Xincrement1}
    X_{n+1} = \begin{cases}
      \tilde{X}_{n+1} & \tilde{X}_{n+1}\in U_{1}\\
      X_{0}^{(\sigma_n)} & \tilde{X}_{n+1} \notin U_{1}
    \end{cases}\\
    \label{e:sigincrement1}
    \sigma_{n+1} = \begin{cases}
      \sigma_n &  \tilde{X}_{n+1} \in U_{1}\\
      \sigma_n+1 &\tilde{X}_{n+1} \notin U_{1}
    \end{cases}
  \end{gather}
\end{subequations}
We interpret $\tilde{X}_{n+1} $ as the proposed move, which is either
accepted or rejected depending on whether or not it will remain in the
trust region.  If it is rejected, the algorithm restarts at
$X_{0}^{(\sigma_n)}\in U_0$. The restart points,
$\{X_{0}^{(\sigma_n)}\}$, may be random, or it may be
that $X_{0}^{(\sigma_n)}= X_0$ is fixed.  The essential property is that the
algorithm will restart in the interior of the trust
region, away from its boundary.   The r.v. $\sigma_n$ counts
the number of times a truncation has occurred. Algorithm
\eqref{e:fixed_trust} can now be expressed as
\begin{equation}
  \label{e:Xfixed_trust}
  \begin{split}
    X_{n+1} &= X_n - a_{n+1} F(X_n, Z_n+1) + P_{n+1}\\
    P_{n+1} & = \{{X_{0}^{(\sigma_n)}
      -\tilde{X}_{n+1}}\}1_{\tilde{X}_{n+1} \notin U_1}.
  \end{split}
\end{equation}

\subsubsection{Expanding Trust Regions}

In the second version of truncated Robbins-Monro, define the sequence of open bounded sets, $U_{n}$ such that:
\begin{gather}
  \label{e:nested}
  U_0\subsetneq U_1\subsetneq U_2\subsetneq\ldots, \quad   \cup_{n=0}^\infty U_n = \calH.
\end{gather}
Again, letting $X_0 \in U_0$, $\sigma_0 = 0$, the algorithm  is
\begin{subequations}
  \label{e:expanding_trust}
  \begin{gather}
    \label{e:proposa2l}
    \tilde{X}_{n+1} = X_n - a_{n+1} F(X_n, Z_{n+1})\\
    \label{e:Xincrement2}
    X_{n+1} = \begin{cases}
      \tilde{X}_{n+1} & \tilde{X}_{n+1} \in U_{\sigma_n}\\
      X_{0}^{(\sigma_n)} & \tilde{X}_{n+1} \notin U_{\sigma_n}
    \end{cases}\\
    \label{e:sigincrement2}
    \sigma_{n+1} = \begin{cases}
      \sigma_n &  \tilde{X}_{n+1} \in U_{\sigma_n}\\
      \sigma_n+1 &\tilde{X}_{n+1}\notin U_{\sigma_n}
    \end{cases}
  \end{gather}
\end{subequations}
A consequence of this formulation is that $X_n \in U_{\sigma_n}$ for all $n$.
As before, the restart points may be random
or fixed, and they are in $U_0$.  This would appear
superior to the fixed trust region algorithm, as it does not require
knowledge of the sets.  However, to guarantee convergence, global (in
$\calH$) assumptions on the regression function are required; see
Assumption \ref{a:convex} below.  \eqref{e:expanding_trust} can
written with $P_{n+1}$ as
\begin{equation}
  \label{e:Xexpanding_trust}
  \begin{split}
    X_{n+1} &= X_n - a_{n+1} F(X_n, Z_n+1) + P_{n+1}\\
    P_{n+1} & = \{{X_{0}^{(\sigma_n)}
      -\tilde{X}_{n+1}}\}1_{\tilde{X}_{n+1} \notin U_{\sigma_n}}
  \end{split}
\end{equation}

\subsection{Outline}

In Section \ref{s:convergence}, we state sufficient assumptions for which we are
able to prove convergence in both the fixed and expanding trust region problems,
and we also establish some preliminary results.  In Section \ref{s:kl}, we focus
on the relative entropy minimization problem, and identify what assumptions must
hold for convergence to be guaranteed.  Examples are then presented in Section
\ref{s:examples}, and we conclude with remarks in Section \ref{s:disc}.

\section{Convergence of Robbins-Monro}
\label{s:convergence}

We first reformulate \eqref{e:Xfixed_trust} and
\eqref{e:Xexpanding_trust} in the more general form
\begin{equation}
  \label{e:Xinc_general}
  X_{n+1} = \underbrace{X_n - a_{n+1} f(X_n) - a_{n+1}\delta M_{n+1}}_{=\tilde{X}_{n+1}} + a_{n+1} P_{n+1},\\
\end{equation}
where $\delta M_{n+1}$, the noise term, is
\begin{equation}
  \label{e:deltaM}
  \begin{split}
  \delta M_{n+1} &= F(X_n, Z_{n+1})   - f(X_n)\\
  &= F(X_n, Z_{n+1})-\E[ F(X_n, Z_{n+1})\mid X_n].
    \end{split}
\end{equation}
A natural filtration for this problem is   $\calF_n = \sigma(X_0, Z_1,
\ldots, Z_n)$.  $X_n$ is $\calF_n$ measurable and the noise term can be expressed
in terms of the filtration as $\delta M_{n+1} = F(X_n, Z_{n+1}) - \E[F(X_n, Z_{n+1})\mid \calF_n]$.

We now state our main assumptions:
\begin{assume}
  \label{a:root}
  $f$ has a zero, $x_\star$. In the case of the fixed trust region problem,
  there exist $R_0<R_1$ such that
  \begin{equation*}
    U_0\subseteq B_{R_0}(x_\star)\subset B_{R_1}(x_\star)\subseteq U_1.
  \end{equation*}
In the case of the expanding trust region problem, the open sets are defined as $U_n = B_{r_n}(0)$ with 
\begin{equation}
\label{e:radii}
    0< r_0< r_1<r_2<\ldots<r_n \to \infty.
\end{equation}
These sets clearly satisfy \eqref{e:nested}.
\end{assume}

\begin{assume}
  \label{a:convex}
  For any $0<a<A$, there exists $\delta >0$:
  \begin{align*}
    \inf_{a\leq\norm{x-x_\star}\leq A} \inner{x-x_\star}{f(x)}\geq \delta.
  \end{align*}
In the case of the fixed truncation, this inequality is restricted to
$x\in U_1$.  This is akin to a convexity condition on a functional $\mathcal{F}$ with $f = D\mathcal{F}$.
\end{assume}

\begin{assume}
  \label{a:bound}
  $x\mapsto \E[\norm{F(x,Z)}^2]$ is bounded on bounded sets, with the
  restriction to $U_1$ in the case of fixed trust
  regions.
\end{assume}

\begin{assume}
  \label{a:sequence}
  $a_n>0$, $\sum a_n = \infty$, and $\sum a_n^2 < \infty$
\end{assume}

\begin{thm}
  \label{t:fixed_trust}
  Under the above assumptions, for the fixed trust region problem,
  $X_n \to x_\star$ a.s. and $\sigma_n$ is a.s. finite.
\end{thm}

\begin{thm}
  \label{t:expanding_trust}
  Under the above assumptions, for the expanding trust region problem,
  $X_n \to x_\star$ a.s. and $\sigma_n$ is a.s. finite.
\end{thm}

Note the distinction between the assumptions in the two algorithms.
In the fixed truncation algorithm, Assumptions \ref{a:convex} and
\ref{a:bound} need only hold in the set $U_1$, while in the expanding
truncation algorithm, they must hold in all of $\calH$.  While this
would seem to be a weaker condition, it requires identification of the
sets $U_0$ and $U_1$ for which the assumptions hold.  Such sets may not
be readily identifiable, as we will see in our examples.


We first need some additional information about $f$ and the noise
sequence $\delta M_n$.
\begin{lem}
  \label{l:boundedness}
  Under Assumption \ref{a:bound}, $f$ is bounded on $U_1$, for the
  fixed trust region problem, and on arbitrary bounded sets, for the
  expanding trust region problem.
\end{lem}

\begin{proof}
  Trivially,
  \begin{equation*}
    \norm{f(x)} = \norm{\E[F(x,Z)]}\leq \E[\norm{F(x,Z)}
    ]\leq \sqrt{\E[\norm{F(x,Z)}^2}],
  \end{equation*}
  and the results follows from the assumption.
\end{proof}

\begin{prop}
  \label{p:mds}
  For the fixed trust region problem, let
  \[
  M_n = \sum_{i=1}^n a_i \delta M_i.
  \]
  Alternatively, in the expanding trust region problem, for $r>0$,
  let
  \[
  M_n = \sum_{i=1}^n a_i \delta M_i 1_{\norm{X_{i-1}-x_\star}\leq r}.
  \]
  Under Assumptions \ref{a:bound} and \ref{a:sequence}, $M_n$ is
  a martingale, converging in $\calH$, a.s.
\end{prop}
\begin{proof}  The following argument holds in both the fixed and expanding trust region problems, with appropriate modifications.  We present the expanding trust region case.  The proof is broken up into \ref{prop21lsast} steps:

\begin{enumerate}[\upshape 1.]

\item Relying on Theorem 6 of \cite{Chatterji:1968un} for Banach space
  valued martingales, it will be sufficient to show that $M_n$
  is a martingale, uniformly bounded in $L^1(\mathbb{P})$. 

  \item  In the case of the expanding truncations,
  \begin{equation*}
    \begin{split}
      \E[\norm{\delta M_i 1_{\norm{X_{i-1}-x_\star}\leq r}}^2] &\leq 2
      \E[\norm{F(X_{i-1},Z_i)1_{\norm{X_{i-1}-x_\star}\leq r}}^2] \\
      &\quad + 2 \E[\norm{f(X_{i-1})1_{\norm{X_{i-1}-x_\star}\leq r}}^2] \\
      &\leq 2\sup_{\norm{x- x_\star}\leq r}
      \E[\norm{F(x,Z)}^2] + 2\sup_{\norm{x- x_\star}\leq
        r}\norm{f(x)}^2
    \end{split}
  \end{equation*}
  Since both of these terms are bounded, independently of $i$, by Assumption
  \ref{a:bound} and Lemma \ref{l:boundedness}, this is finite.

  \item Next, since $\{\delta M_i 1_{\|X_{i-1}-x_\star\|\leq r}\}$ is a
  martingale difference sequence, we can use the above estimate to
  obtain the uniform $L^2(\mathbb{P})$ bound,
  \begin{equation*}
    \begin{split}
      \E[\norm{M_n}^2]& = \sum_{i=1}^na_i^2 \E[\norm{{\delta M_i
          1_{\norm{X_{i-1}-x_\star}\leq r}}}^2] \\
      &\leq\sup_{i}\E[\norm{{\delta M_i 1_{\norm{X_{i-1}-x_\star}\leq
            r}}}^2]\sum_{i=1}^\infty a_i^2 < \infty
    \end{split}
  \end{equation*}
  Uniform boundedness in $L^2$, gives boundedness in
  $L^1$, and this implies a.s. convergence in $\calH$.
  \label{prop21lsast}

\end{enumerate}

\end{proof}

\subsection{Finite Truncations} In this section we prove results showing that
only finitely many truncations will occur, in either the fixed or expanding
trust region case.  Recall that when a truncation occurs, the equivalent conditions hold: $P_{n+1}\neq 0$; $\sigma_{n+1} = \sigma_n+1$; and $\tilde{X}_{n+1}\notin U_1$ in the fixed trust region algorithm, while  $\tilde{X}_{n+1}\notin U_{\sigma_n}$ in the expanding trust region case.

\begin{lem}
  \label{l:fixed_finite_truncations}

  In the fixed trust region
algorithm, if Assumptions \ref{a:root}, \ref{a:convex}, \ref{a:bound}, and
\ref{a:sequence} hold, then the number of truncations is a.s. finite; a.s.,
there exists $N$, such that for all $n\geq N$, $\sigma_n = \sigma_N$.

\end{lem}

\begin{proof}

We break the proof up into \ref{lemma22last} steps:

\begin{enumerate}[\upshape 1.]
  \item Pick  $\rho$ and $\rho'$ such that
  \begin{equation}
  \label{e:Rrho}
  R_0<\rho' < \rho < R_1
  \end{equation}
  Let $\bar f = \sup\|f(x)\|$, with the supremum over $U_1$; this
  bound exists by Lemma \ref{l:boundedness}. Under Assumption \ref{a:convex}, there exists $\delta>0$ such that
  \begin{equation}
    \label{e:convexR0}
    \inf_{R_0/2 \leq \|x-x_\star\|\leq R_1}\inner{x-x_\star}{f(x)}= \delta.
  \end{equation}
 Having fixed $\rho$, $\rho'$, $\bar f$, and $\delta$,
  take $\eps>0$  such that:
\begin{equation}
  \label{e:lemcond}
  \eps < \min\set{\rho'-R_0, \frac{R_1-\rho'}{2 + \bar f}, \frac{\rho' - R_0}{\bar f}, \frac{R_0}{2}, \frac{\delta}{2\bar{f}}, \frac{\delta}{\bar{f}^2}, {\rho - \rho'}}.
\end{equation}
   Having fixed such an $\eps$, by the assumptions of this
  lemma and Proposition \ref{p:mds}, a.s., there exists $n_{\eps}$ such that for any
  $n,m\geq n_{\eps}$, both
  \begin{equation}
  \label{e:seqcond}
    \norm{\sum_{k=n}^m a_k \delta M_k}\leq \eps, \quad a_n \leq \eps.
  \end{equation}

  \item  Define the auxiliary sequence
  \begin{equation}
    \label{e:Xprimeseq1}
    X_n' = X_n - \sum_{k=n+1}^\infty a_k \delta M_k.
  \end{equation}
  Using \eqref{e:Xinc_general}, we can then write
  \begin{equation}
    \label{e:Xprimeinc1}
    X_{n+1}' = X_{n}' - a_{n+1}f(X_n) +a_{n+1}P_{n+1}.
  \end{equation}
  By \eqref{e:seqcond}, for any $n \geq n_\eps$,
  \begin{equation}
  \label{e:seqdiff}
      \|X_n' - X_n\|\leq \eps
  \end{equation}

\item We will show  $X_n' \in B_{\rho'}(x_\star)$ for all $n$ large enough.  The significance of this is that if $n\geq n_\eps$, and $X_n'\in B_{\rho'}(x_\star)$, then no truncation occurs.  Indeed, using \eqref{e:lemcond}
     \begin{equation}
     \label{e:XauxP0}
    \begin{split}
      \|\tilde{X}_{n+1} -x_\star\| &\leq \|X_{n}' -x_\star\|
      +\|X_n - X_n'\| + a_{n+1}\bar f  + \|a_{n+1}\delta M_{n+1}\|\\
      &< \rho' + \eps + \eps \bar f + \eps<R_1, \Rightarrow \tilde{X}_{n+1}\in U_1.
    \end{split}
  \end{equation}
  Consequently, $P_{n+1}=0$, $X_{n+1} = \tilde{X}_{n+1}$, and $\sigma_{n+1} = \sigma_n$.  Thus, establishing $X_n'\in B_{\rho'}(x_\star)$ will yield the result.

\item  Let
\begin{equation}
\label{e:Ntruncation}
  N = \inf\{n\geq n_\eps\mid \tilde{X}_{n+1}\not\in U_1\}+1
\end{equation}
This corresponds to the the first truncation after $n_\eps$.  If the above set
is empty, for that realization, no truncations occur after $n_\eps$, and we are
done.  In such a case, we may take $N = n_\eps$ in the statement of the lemma.

\item We now prove by induction that in the case that \eqref{e:Ntruncation} is finite, $X_n'\in B_{\rho'}(x_\star)$ for all $n\geq N$.  First, note that $X_N \in B_{R_0}(x_\star)\subset B_{\rho}(x_\star)$.    By \eqref{e:lemcond} and \eqref{e:seqdiff},
  \[
  \|X_{N}'-x_\star\| \leq \|X_{N}-x_\star\| + \|X_{N}'- X_{N}\|< R_0
  + \eps < \rho', \Rightarrow X_N' \in B_{\rho'}(x_\star).
  \]

    Next, assume
  $X_{N}', X_{N+1}', \ldots, X_n'$ are all in $B_{\rho'}(x_\star)$.  Using \eqref{e:XauxP0}, we have that $P_{N+1}=\ldots=P_{n+1}=0$ and $\sigma_{N}=\ldots= \sigma_{n}$.  Therefore,
  \begin{equation}
    \label{e:fixedtrustlem1}
    \begin{split}
      \norm{X_{n+1}' - x_\star}^2 &= \norm{X'_n - x_\star}^2 - 2a_{n+1}
      \inner{X'_{n}-x_\star}{f(X_n)} \\
      &\quad + a_{n+1}^2 \norm{f(X_n)}^2\\
      &\leq \norm{X'_n - x_\star}^2 - 2a_{n+1}
      \inner{X'_{n}-x_\star}{f(X_n)} + a_{n+1}\eps \bar f^2
    \end{split}
  \end{equation}
  We now consider two cases of \eqref{e:fixedtrustlem1} to conclude  $\|X_{n+1}' - x_\star\|< \rho'$.

  \item  In the first case, $\|X'_n - x_\star\|\leq R_0$.  By Cauchy-Schwarz and \eqref{e:lemcond}
  \[
  \|X_{n+1}' - x_\star\|^2< R_0^2 + 2 \eps R_0 \bar f + \eps^2\bar{f}^2 = (R_0 + \eps \bar f)^2< (\rho')^2.
  \]

 In the second case, $R_0<\|X'_n - x_\star\|< \rho'$. Dissecting the inner product term in \eqref{e:fixedtrustlem1} and
  using Assumption \ref{a:convex} and \eqref{e:seqdiff},
  \begin{equation}
  \label{e:ip_split}
    \begin{split}
      \inner{X'_{n}-x_\star}{f(X_n) } &= \inner{X_{n}-x_\star}{f(X_n)
      } +
      \inner{X'_{n}-X_n}{f(X_n) }\\
      &\geq\inner{X_{n}-x_\star}{f(X_n)
      }-\bar f\eps
    \end{split}
  \end{equation}
  Conditions \eqref{e:lemcond} and \eqref{e:seqdiff} yield the following upper and lower bounds:
  \begin{align*}
  \|X_n-x_\star\|&\geq\|X_n'-x_\star\| -\|X_n'-X_n\| \geq  R_0 -\eps> \tfrac{1}{2}R_0, \\
  \|X_n-x_\star\|&\leq \|X_n'-x_\star\| +\|X_n'-X_n\| \leq  \rho' +\eps< \rho< R_1.
  \end{align*}
  Therefore, \eqref{e:convexR0} applies and $\inner{X_{n}-x_\star}{f(X_n)}\geq \delta$. Using this in \eqref{e:ip_split}, and condition \eqref{e:lemcond},
   \begin{equation*}
       \inner{X'_{n}-x_\star}{f(X_n) }\geq \delta - \bar f \eps >\tfrac{1}{2}\delta.
   \end{equation*}
  Substituting this last estimate  back into \eqref{e:fixedtrustlem1}, and using \eqref{e:lemcond},
  \[
  \|X_{n+1}' - x_\star\|^2< (\rho')^2 -  a_{n+1}(\delta - \eps\bar f^2)<(\rho')^2.
  \]
  This completes the inductive step.

  \item  Since the auxiliary sequence remains in $B_{\rho'}(x_\star)$ for all $n\geq N>n_\eps$, \eqref{e:XauxP0} ensures $\tilde X_{n+1}\in B_{R_1}(x_\star)$, $P_{n+1}=0$, and $\sigma_{n+1}= \sigma_N$, a.s.
  
  \label{lemma22last}

\end{enumerate}

\end{proof}

To obtain a similar result for the expanding trust region problem, we
first relate the finiteness of the number of truncations with the
sequence persisting in a bounded set.
\begin{lem}
  \label{l:bound_finite}
  In the expanding trust region algorithm, if Assumptions \ref{a:root},
\ref{a:bound}, and \ref{a:sequence} hold, then the sequence remains in a set of
the form $B_{R}(0)$ for some $R>0$ if and only if the number of truncations is
finite, a.s.
\end{lem}

\begin{proof}
We break this proof into \ref{lemma23last} steps:
\begin{enumerate}[\upshape 1.]

\item If the number of truncations is finite, then there exists $N$ such that for all $n\geq N$, $\sigma_n = \sigma_{N}$.  Consequently, the proposed moves are always
accepted, and $X_{n} \in U_{\sigma_{n}}= U_{\sigma_N}$ for all $n\geq N$.  Since $X_n \in U_{\sigma_n}\subset U_{\sigma_N}$ for $n< N$, $X_n \in U_{\sigma_N}$ for all $n$.  By Assumption \ref{a:bound}, $B_R(0)=B_{r_{\sigma_N}}(0)=U_{\sigma_N}$ is the desired set.

\item  For the other direction, assume that  there exists $R>0$ such that $X_n \in
B_R(0)$ for all $n$.  Since the $r_n$ in \eqref{e:radii}
tend to infinity, there exists $N_1$, such that $R < R +1 < r_{N_1}$.  Hence,
for all $n\geq N_1$,
  \begin{equation}
  \label{e:Brsubset}
   B_{R}(0)\subset B_{R+1}(0) \subset U_{n}
  \end{equation}
  Let $\bar f = \sup \|f(x)\|$, with the supremum over $B_{R}(0)$.  Let $\tilde{R}$ be sufficiently large such that $B_{R+1}(0)\subset B_{\tilde R}(x_\star)$.  Lastly, using Proposition \ref{p:mds} and Assumption \ref{a:sequence}, a.s., there exists $N_2$, such that for all $n\geq N_2$
  \begin{equation}
  \label{e:seqcond2}
     \| a_{n} \delta M_{n}1_{\|X_n - x_\star\|\leq \tilde R}\|<\frac{1}{2}, \quad a_{n}< \frac{1}{2(1+ \bar f)}
 \end{equation}
 Since  $X_n \in B_{R}(0) \subset B_{\tilde R}(x_\star)$, the indicator function in \eqref{e:seqcond2} is always one, and $\| a_{n} \delta M_{n}\|<1/2$.

 \item Next, let
 \begin{equation}
 \label{e:lemma23set}
     N = \inf \{ n\geq 0\mid \sigma_n\geq \max\{N_1, N_2\}\}
 \end{equation}
 If the above set is empty, then $\sigma_n< \max\{N_1, N_2\}$ for all $n$, and the number of truncations is a.s. finite.  In this case, the proof is complete.

 \item If the set in \eqref{e:lemma23set} is not empy, then $N < \infty$.  Take $n\geq N$.   As $X_n \in B_R(0)$, and since $n\geq \sigma_n\geq \max\{N_1, N_2\}$, \eqref{e:seqcond2} applies.  Therefore,
 \begin{equation}
 \begin{split}
     \|\tilde{X}_{n+1}\| &\leq \|X_n\| + \| \tilde{X}_{n+1} - X_n\|\\
     &\leq \|X_n\|  + a_{n+1}\|f(X_n)\| + \| a_{n+1} \delta M_{n+1}\|\\
     &< R + \tfrac{1}{2}+\tfrac{1}{2}< R +1.
     \end{split}
 \end{equation}
 Thus, $\tilde{X}_{n+1}\in B_{R+1}(0)\subset U_{N_1}$, $\sigma_n \geq N_1$,  and $U_{N_1}\subset U_{\sigma_n}$.  Therefore, $\tilde{X}_{n+1}\in U_{\sigma_n}$.  No truncation occurs, and $\sigma_n = \sigma_{n+1}$. Since this holds for all $n\geq N$, $\sigma_n = \sigma_N$, and the number of truncations is a.s. finite.
 
 \label{lemma23last}
\end{enumerate}

\end{proof}

Next, we establish that, subject to an additional assumption, the
sequence remains in a bounded set; the finiteness of the truncations
is then a corollary.
\begin{lem}
  \label{l:bounded1}
  In the expanding trust region algorithm, if Assumptions \ref{a:root}, \ref{a:convex}, \ref{a:bound}, and \ref{a:sequence}
  hold, and for any $r>0$, there a.s. exists $N<\infty$,
  such that for all $n \geq N$,
  \begin{equation*}
    P_{n+1} 1_{\|X_n- x_\star\|\leq r}=0,
  \end{equation*}
  then $\{X_n\}$ remains in a bounded open set, a.s.
\end{lem}

\begin{proof}
We  break this proof into \ref{lemma24last} steps:
\begin{enumerate}[\upshape 1.]

\item We begin by setting some constants for the rest of the proof.  Fix $R>0$ sufficiently large such that $B_{R}(x_\star)\supset U_0$.   Next, let
  $\bar f = \sup {\|f(x)\|}$ with the supremum taken over $B_{R+2}(x_\star)$.  Assumption \ref{a:convex} ensures there exists $\delta>0$ such that
  \begin{equation}
    \label{e:convexR1}
    \inf_{R/2 \leq \|x-x_\star\|\leq R+2}\inner{x-x_\star}{f(x)}= \delta.
  \end{equation}
  Having fixed $R$, $\bar f$, and $\delta$, take $\eps >0$
  such that:
  \begin{equation}
    \label{e:lem24cond}
    \eps < \min\set{1,\frac{1}{\bar{f}}, \frac{\delta}{2\bar{f}}, \frac{\delta}{\bar{f}^2}, \frac{R}{2}}.
  \end{equation}
  By the assumptions of this lemma and Proposition \ref{p:mds} there exists, a.s.,
$n_\eps\geq N$ such that for all $n\geq n_\eps$,
  \begin{subequations}
    \label{e:seqcond3}
    \begin{gather}
    \norm{\sum_{i=n+1}^\infty a_i \delta M_i 1_{\norm{X_{i-1} -
          x_\star}\leq
        R+2}}\leq \eps,\\
         P_{n+1} 1_{\norm{X_{n}- x_\star}\leq R+2}=0,\\
         a_{n+1}\leq \eps
    \end{gather}
  \end{subequations}

\item Define the modified sequence for $n \geq n_\eps$ as
  \begin{equation}
    \label{e:modseq2}
    X_n' = X_n - \sum_{k=n+1}^\infty a_k \delta M_k 1_{\norm{X_{k-1} -
        x_\star}\leq R+2}, \Rightarrow \|X_n' - X_n\|\leq \eps.
  \end{equation}
  Using \eqref{e:Xinc_general}, we have the iteration
  \begin{equation}
    \label{e:moditer2}
    X_{n+1}' = X_n' - a_{n+1} \delta M_{n+1} 1_{\norm{X_n - x_\star}> R+2}
    - a_{n+1}f(X_n) +a_{n+1} P_{n+1}.
  \end{equation}

\item  Let
\begin{equation}
\label{e:N2}
    N = \inf\{n\geq n_\eps\mid \sigma_{n+1}\neq \sigma_n \}+1,
\end{equation}
the first time  after $n_\eps$ that a truncation occurs.  

If the above set is empty, no truncations occur after $n_\eps$.  In this case, $\sigma_n =
\sigma_{n_\eps}\leq n_\eps<\infty$  for all $n \geq n_\eps$.  Therefore, for all
$n\geq n_\eps$, $X_n \in U_{\sigma_n}\subset U_{\sigma_{n_\eps}}$.    Since $U_{\sigma_n}\subset
U_{\sigma_{n_\eps}}$ for all $n < n_\eps$ too, the proof is complete in this case.

\item Now assume that $N<\infty$.  We will show that  $\{X_n'\}$ remains in  $B_{R+1}(x_\star)$ for all $n\geq N$.  Were this to hold, then for $n\geq N$,
\begin{equation}
\label{e:seqR2}
\begin{split}
    \|X_n - x_\star\| &\leq \|X_n' - x_\star\| + \norm{\sum_{i=n+1}^\infty a_i \delta M_i 1_{\norm{X_{i-1} -
        x_\star}\leq R+2}}\\
        &< R+1 + \eps < R+2,
\end{split}
\end{equation}
having used \eqref{e:seqcond3} and \eqref{e:modseq2}.  For $n < N$, $X_n \in U_{\sigma_n}\subset U_{\sigma_N} = B_{r_N}(0)$.  Therefore, for all $n$,  $X_n \in B_{\tilde{R}}(0)$ where $\tilde{R} = \max\{ r_N, \|x_\star\| + R+2\}$.

\item   We prove $X_n' \in B_{R+1}(x_\star)$ by induction.  First,
since $\eps <1$ and $X_N \in U_{0}\subset B_R(x_\star)$,
\begin{equation*}
    \|X_N'-x_\star\|\leq \|X_N'-X_N\|+ \|X_N-x_\star\|< \eps + R< R+1.
\end{equation*}

Next, assume that $X_N',X_{N+1}',\ldots, X_n'$ are all in $B_{R+1}(x_\star)$. By \eqref{e:seqR2}, $X_n\in B_{R+2}(x_\star)$. Since $P_{n+1} 1_{\norm{X_{n}- x_\star}\leq R+2}=0$, we conclude  $P_{n+1}=0$.  The modified iteration \eqref{e:moditer2} simplifies to
  have
  \[
  X_{n+1}' = X_n' -a_{n+1}f(X_n),
  \]
  and
  \begin{equation}
  \label{e:lemming_bound1}
    \begin{split}
      \|X_{n+1}' - x_\star\|^2 &= \|X'_n - x_\star\|^2 - 2a_{n+1}
      \inner{X'_{n}-x_\star}{f(X_n) }\\
      &\quad + a_{n+1}^2 \|f(X_n)\|^2\\
      &< \|X'_n - x_\star\|^2 - 2a_{n+1}
      \inner{X'_{n}-x_\star}{f(X_n) }\\
      &\quad + a_{n+1} \eps \bar{f}^2.
    \end{split}
  \end{equation}

\item  We now consider two cases of \eqref{e:lemming_bound1}.  First, assume $\|X'_{n} - x_\star\|\leq R$.  Then
   \eqref{e:lemming_bound1} can immediately be bounded as
  \begin{equation*}
  \begin{split}
    \|X_{n+1}' - x_\star\|^2&< R^2 + 2 \eps R \bar{f} +\eps^2\bar{f}^2 = (R+\eps \bar{f})^2< (R+1)^2,
    \end{split}
  \end{equation*}
  where we have used condition \eqref{e:lem24cond} in the last inequality.

  \item  Now  consider the case $R<\|X'_{n} - x_\star\|< R+1$.  Using \eqref{e:lem24cond}, the inner product in \eqref{e:lemming_bound1} can first be bounded from below:
  \begin{equation*}
  \begin{split}
      \inner{X'_{n}-x_\star}{f(X_n) } &= \inner{X_{n}-x_\star}{f(X_n)
      } +
      \inner{X'_{n}-X_n}{f(X_n) }\\
      &\geq \inner{X_{n}-x_\star}{f(X_n)
      } - \eps \bar{f}> \inner{X_{n}-x_\star}{f(X_n)
      } - \tfrac{1}{2}\delta.
  \end{split}
\end{equation*}
 Next, using \eqref{e:lem24cond}
  \begin{equation*}
      \|X_{n} - x_\star\|\geq \|X_{n}'-x_\star\| - \|X_{n} -X'_{n}\|> R-\eps>R - \tfrac{1}{2}R=\tfrac{1}{2}R
  \end{equation*}
  Therefore, $\tfrac{1}{2}R<\|X_n -x_\star\| < R+2$, so \eqref{e:convexR1} ensures $  \inner{X_{n}-x_\star}{f(X_n)}\geq\delta$ and
  \begin{equation*}
       \inner{X'_{n}-x_\star}{f(X_n) }>\delta-\tfrac{1}{2}\delta=\tfrac{1}{2}\delta.
  \end{equation*}

    Returning to \eqref{e:lemming_bound1}, by \eqref{e:lem24cond},
  \begin{equation*}
    \|X_{n+1}' - x_\star\|^2\leq (R+1)^2 - a_{n+1}(\delta - \eps \bar{f}^2)< (R+1)^2.
  \end{equation*}
  This completes the proof of the inductive step in this second case, completing the proof.
  \label{lemma24last}
\end{enumerate}
\end{proof}

\begin{cor}
  \label{c:finite_truncations}
  For the expanding trust region algorithm, if Assumptions \ref{a:root},
\ref{a:convex}, \ref{a:bound}, and \ref{a:sequence} hold, then the number of
truncations is a.s. finite.
\end{cor}

\begin{proof}
  The proof is by contradiction.  We break the proof into \ref{cor21last} steps:
\begin{enumerate}[\upshape 1.]

  \item Assuming that there are infinitely many truncations, Lemma
\ref{l:bound_finite} implies that the sequence cannot remain in a bounded set.
Then, continuing to assume that  Assumptions \ref{a:root}, \ref{a:convex},
\ref{a:bound}, and \ref{a:sequence} hold, the only way for the conclusion of
Lemma \ref{l:bounded1} to fail is if the assumption on $P_{n+1} 1_{\|X_n -
x_\star\|\leq r }$ is false.    Therefore, there exists $r>0$ and a set of
positive measure on which a subsequence, $P_{n_k+1} 1_{\|X_{n_k} - x_\star\|\leq
r }\neq 0$.   Hence $X_{n_k}\in B_{r}(x_\star)$, and $P_{n_k+1}\neq 0$.  So
truncations occur at these indices, and $\tilde{X}_{n_k+1} \not\in
U_{\sigma_{n_k+1}}$.

    \item Let $\bar{f} = \sup \|f(x)\|$ with the  supremum over the set
    $B_r(x_\star)$ and let $\eps>0$ satisfy
    \begin{equation}
    \label{e:cor_cond1}
        \eps < (\bar{f}+1)^{-1}.
    \end{equation}
    By our assumptions of the lemma and Proposition \ref{p:mds}, there exists $n_\eps$ such that for all
    $n\geq n_\eps$
    \begin{equation}
        \|a_{n+1}\delta M_{n+1}1_{\|X_n - x_\star\|\leq r}\|\leq \eps, \quad a_{n+1}\leq \eps
    \end{equation}
    Along the subsequence, for all  $n_k\geq n_\eps$,
    \begin{equation}
        \|a_{n_k+1}\delta M_{n_k+1}1_{\|X_{n_k} - x_\star\|\leq r}\|=\|a_{n_k+1}\delta M_{n_k+1}\|\leq \eps.
    \end{equation}

\item Furthermore, for $n_k \geq n_\eps$:
\begin{equation}
\label{e:corbound}
\begin{split}
    \|\tilde{X}_{n_k+1} - x_\star\| &\leq \| X_{n_k} - x_\star\| + a_{n_k+1}\|f(X_{n_k})\| + \|a_{n_k+1}\delta M_{n_k+1}\|\\
    &< r + \eps \bar{f} + \eps < r+1, \Rightarrow \tilde{X}_{n_k+1} \in B_{r+1}(x_\star),
\end{split}
\end{equation}
where \eqref{e:cor_cond1} has been used in the last inequality.

\item  By the definition of the $U_n$, there exists an index $M$ such that $U_{M}\supset B_{r+1}(x_\star)$.  Let
\begin{equation}
    N = \inf\{n\geq n_\eps\mid \sigma_n \geq M\}.
\end{equation}
This set is nonempty and $N<\infty$, since we have assumed there are infinitely many truncations.  Let $n_k \geq N$.  Then $\sigma_{n_k}\geq M$ and $U_{\sigma_{n_k}}\supset B_{r+1}(x_\star)$.  But \eqref{e:corbound} then implies that $\tilde{X}_{n_k+1} \in U_{\sigma_{n_k}}$, and no truncation will occur; $P_{n_k+1}=0$, providing the desired the contradiction.

\label{cor21last}

\end{enumerate}
\end{proof}

\subsection{Proof of Convergence}
Using the above results, we are able to prove Theorems
\ref{t:fixed_trust} and \ref{t:expanding_trust}.  Since the proofs are
quite similar, we present the more complicated expanding trust
region case.

\begin{proof}  We split this proof into \ref{thmlast} steps:

\begin{enumerate}[\upshape 1.]

\item First, by Corollary \ref{c:finite_truncations}, only finitely many truncations
occur.   By Lemma \ref{l:bound_finite}, there exists $R>0$ such that  $X_n\in
B_R(0)$ for all $n$.  Consequently, there is an $r$ such that
$X_n\in B_r(x_\star)$ for all $n$.

\item Next, we fix  constants.  Let $\bar f = \sup \|f(x)\|$ with the supremum
taken over $B_r(x_\star)$.  Fix $\eta \in (0, 2R)$, and use Assumption
\ref{a:convex} to determine $\delta>0$ such that
   \begin{equation}
    \label{e:convexeta}
         \inf_{\eta/2  \leq \|x-x_\star\|\leq r}\inner{x-x_\star}{f(x)}= \delta
  \end{equation}
    Take $\eps >0$  such that:
    \begin{equation}
      \label{e:thmcond}
      \eps < \min\set{1,\frac{\eta}{2}, \frac{\delta}{2\bar{f}}, \frac{\delta}{2\bar{f}^2}}
    \end{equation}
    Having set $\eps$, we again appeal to Assumption \ref{a:sequence} and Proposition \ref{p:mds} to find $n_\eps$ such that for all $n\geq n_\eps$:
    \begin{equation}
    \label{e:seqcondthm}
  \norm{\sum_{i=n+1}^\infty a_i \delta M_i 1_{\|X_{i-1} -x_\star\|\leq
      r}} =\norm{\sum_{i=n+1}^\infty a_i \delta M_i } \leq \eps,
  \quad a_{n+1}\leq \eps
  \end{equation}

  \item  Define  the auxiliary sequence, 
  \begin{equation}
  \label{e:thmauxseq}
  X_n' = X_n - \sum_{i=n+1}^\infty a_i \delta M_i 1_{\|X_{i-1} -
    x_\star \|\leq r} = X_n - \sum_{i=n+1}^\infty a_i \delta M_i.
  \end{equation}
  Since there are only finitely many truncations, there exists $N\geq n_\eps$, such that for all $n\geq N$, $P_{n+1}=0$, as the truncations have ceased.  Consequently, for $n\geq N$,
  \begin{equation}
      \label{e:iterthm}
      X'_{n+1} = X'_n - a_{n+1}f(X_n)
  \end{equation}
  By \eqref{e:seqcondthm} and \eqref{e:thmauxseq}, for $n\geq N$, $\|X_n - X_n'\|\leq \eps$.  Since $\eps>0$ may be arbitrarily small, it will be sufficient to prove  $X_n'\to x_\star$.

  \item To obtain convergence of $X_n'$, we first examine $\|X_n'-x_\star\|$.  For $n\geq N$,
    \begin{equation}
    \label{e:thm1}
    \begin{split}
      \|X_{n+1}' - x_\star\|^2 &\leq \|X_{n}' -
        x_\star\|^2-2a_{n+1}\inner{X_n' - x_\star}{f(X_n)} + a_{n+1}\eps
      \bar f^2,
    \end{split}
  \end{equation}
  Now consider two cases of this expression.  First, assume  $\|X_n' - x_\star\|\leq\eta$.  In this case, using \eqref{e:thmcond},
  \begin{equation}
  \label{e:thm_regime1}
  \begin{split}
       -2a_{n+1}\inner{X_n' - x_\star}{f(X_n)} + a_{n+1}\eps
      \bar f^2&\leq a_{n+1}(\eps \bar f + \eps \bar f^2) \\
      &< a_{n+1}(\bar f + \bar f^2)= a_{n+1} B.
  \end{split}
  \end{equation}
  where $B>0$ is a constant depending only on   $\bar f$.
  For $\|X_n' - x_\star\|>\eta$, using \eqref{e:thmcond}
  \begin{equation}
  \label{e:thmip}
  \begin{split}
    \inner{X'_{n}-x_\star}{f(X_n) } &= \inner{X_{n}-x_\star}{f(X_n) } +  \inner{X'_{n}-X_n}{f(X_n) }\\
    &\geq \inner{X_{n}-x_\star}{f(X_n) } - \eps \bar{f}\\
    &> \inner{X_{n}-x_\star}{f(X_n) } -\tfrac{1}{2}\delta.
    \end{split}
  \end{equation}
  By \eqref{e:thmcond},
  \begin{equation*}
      \|X_{n}-x_\star\|\geq \|X'_{n}-x_\star\| -\| X_{n}-X_n'\|> \eta - \eps> \tfrac{1}{2}\eta
  \end{equation*}
  Since $\|X_n - x_\star\|< r$ too, \eqref{e:convexeta} and \eqref{e:thmip} yield the estimate
  \begin{equation*}
      \inner{X'_{n}-x_\star}{f(X_n) }>\delta - \eps \bar{f}>\tfrac{1}{2}\delta
  \end{equation*}
  Thus, in this regime, using \eqref{e:thmcond},
  \begin{equation}
  \label{e:thm_regime2}
  \begin{split}
      -2a_{n+1}\inner{X_n' - x_\star}{f(X_n)} + a_{n+1}\eps
      \bar f^2&\leq -a_{n+1}(\delta - \eps \bar f^2)\\
      &< -\tfrac{1}{2}\delta a_n+1= - A a_{n+1}
  \end{split}
  \end{equation}
  where $A>0$ is a constant depending only on $\delta$.

  Combining estimates \eqref{e:thm_regime1} and \eqref{e:thm_regime2}, we can write for $n\geq N$
  \begin{equation}
  \label{e:thm2}
       \|X_{n+1}' - x_\star\|^2 < \|X_{n}' -
        x_\star\|^2 - a_{n+1} A 1_{\|X_n' - x_\star\|>\eta}  + a_{n+1} B 1_{\|X_n' - x_\star\|\leq \eta}.
  \end{equation}

  \item We now show that $\|X_n' - x_\star\|\leq \eta$ i.o.  The argument is by contradiction.  Let $M\geq N$ be such that for all $n\geq M$, $\|X_n' - x_\star\|> \eta$.  For such $n$,
  \begin{equation}
  \begin{split}
      \eta^2<\|X_{n+1}'-x_\star\|^2  &< \|X_{n}'-x_\star\|^2 - a_{n+1}A \\
      &< \|X_{n-1}'-x_\star\|^2 - a_{n+1}A - a_n A\\
      &<\ldots < \|X_M' - x_\star\|^2 - A \sum_{i=M}^n a_{i+1}.
  \end{split}
  \end{equation}
  Using Assumption \ref{a:sequence} and taking $n\to \infty$, we obtain a contradiction.

  \item Finally, we prove convergence of $X_n'\to x_\star$. Since $X_n'\in B_{\eta}(x_\star)$ i.o., let
  \begin{equation}
      N'= \inf \{n\geq  N\mid \|X_n' - x_\star\|<\eta \}.
  \end{equation}
  For $n\geq N'$, we can then define
  \begin{equation}
    \varphi(n) = \max\set{p\leq n\mid \norm{X_p' - x_\star}< \eta }.
  \end{equation}
  For all such $n$, $\varphi(n)\leq n$, and $X_{\varphi(n)}\in B_\eta(x_\star)$.

  We claim that for  $n\geq N'$,
  \[
  \|X_{n+1}' - x_\star\|^2 < \|X_{\varphi(n)}' -x_\star\|^2 + B
  a_{\varphi(n)+1}< \eta^2 + Ba_{\varphi(n)+1}.
  \]
  First, if $n = \varphi(n)$, this trivially holds in \eqref{e:thm2}.   Suppose now that $n > \varphi(n)$.  Then for
  $i = \varphi(n)+1, \varphi(n)+2,\ldots n$,
  $\|X_{i}' - x_\star\|>\eta$.  Consequently,
  \begin{equation*}
    \begin{split}
      \|X_{n+1}' - x_\star\|^2&< \|X_{n}' - x_\star\|^2< \|X_{n-1}' - x_\star\|^2< \ldots\\
      &< \|X_{\varphi(n)+1}' - x_\star\|^2< \|X_{\varphi(n)}' -
      x_\star\|^2 + B a_{\varphi(n)+1}\\
      & < \eta^2 + B a_{\varphi(n)+1}
    \end{split}
  \end{equation*}
  As $\varphi(n)\to \infty$,
  \begin{equation*}
    \limsup_{n\to \infty }\|X_{n+1}' - x_\star\|^2\leq \eta^2
  \end{equation*}
  Since $\eta$ may be arbitrarily small, we conclude that
  \[
  \limsup_{n\to \infty }\|X_{n+1}' - x_\star\| = \lim_{n\to \infty }\|X_{n+1}' - x_\star\|=0,
  \]
  completing the proof.
  
  \label{thmlast}
\end{enumerate}

\end{proof}

\section{Minimization of Relative Entropy}
\label{s:kl}
Recall from the introduction that our distribution of interest, $\mu$,
is posed on the Borel subsets of Hilbert space $\calH$. We assume that
$\mu \ll \mu_0$, where $\mu_0 = N(m_0, C_0)$ is some reference
Gaussian.  Thus, we write
\begin{equation}
  \label{e:rnd0}
  \frac{d\mu}{d\mu_0} = \frac{1}{Z_\mu}\exp\set{-\Phi_\mu(u)},
\end{equation}
where $\Phi_\nu: X\to \R$, $X$ a Banach space, a subspace of $\mathcal{H}$, of full measure with
respect to $\mu_0$, a Gaussian on $\mathcal{H}$, assumed to be continuous.
$Z_\mu= \E^{\mu_0}[\exp\set{-\Phi(u)}]\in (0,\infty)$ is the
partition function ensuring we have a probability measure.

Let $\nu = N(m, C)$, be another Gaussian, equivalent to $\mu_0$, such
that we can write
\begin{equation}
  \label{e:rnd1}
  \frac{d\nu}{d\mu_0} = \frac{1}{Z_\nu}\exp\set{-\Phi_\nu(v)},
\end{equation}
Assuming that $\nu \ll \mu$, we can write
\begin{equation}
  \label{e:kl1}
  \calR(\nu||\mu) = \E^{\nu}[\Phi_\mu(u) - \Phi_\nu(u)] + \log(Z_\mu) - \log(Z_\nu)
\end{equation}
The assumption that $\nu \ll \mu$ implies that $\nu$ and
$\mu$ are equivalent measures.  As was proven in \cite{Pinski:2015jn},
if $\mathcal{A}$ is a set of Gaussian measures, closed under weak
convergence, such that at least one element of $\mathcal{A}$ is
absolutely continuous with respect to $\mu$, then any minimizing
sequence over $\mathcal{A}$ will have a weak subsequential limit.

If we assume, for this work, that $C=C_0$, then, by the Cameron-Martin
formula (see \cite{Da-Prato:2006aa}),
\begin{equation}
  \label{e:cm_potential}
  \Phi_\nu(u) = -\inner{u-m}{m -m_0}_{\calH^1} - \frac{1}{2}\norm{m -
    m_0}_{\calH^1}^2, \quad Z_{\nu} = 1.
\end{equation}
Here, $\inner{\bullet}{\bullet}_{\calH^1}$ and $\|\bullet\|_{\calH^1}$ are
the inner product and norms of the Cameron-Martin Hilbert space, denoted $\calH^1$,
\begin{equation}
  \label{e:cm_ip}
  \inner{f}{g}_{\calH^1} = \inner{C_0^{-1/2} f}{C_0^{-1/2} g}, \quad
  \norm{f}_{\calH^1}^2 = \inner{f}{f}_{\calH^1}^2.
\end{equation}
Convergence to the minimizer will be established in $\calH^1$, and $\calH^1$ will be the relevant Hilbert space in our application of Theorems \ref{t:fixed_trust} and \ref{t:expanding_trust} to this problem.

Letting $\nu_0 = N(0,C_0)$ and $v\sim \nu_0$, we can then rewrite \eqref{e:kl1} as
\begin{equation}
  \label{e:kl2}
  J(m) \equiv\calR(\nu||\mu) = \E^{\nu_0}[\Phi_{\mu}(v + m)] + \frac{1}{2}\norm{m -
    m_0}_{\calH^1}^2  + \log(Z_\mu)
\end{equation}
The Euler-Lagrange equation associated with \eqref{e:kl2}, and the
second variation, are:
\begin{align}
  \label{e:el1}
  J'(m) &= \E^{\nu_0}[\Phi'_\mu(v+m)] + C_0^{-1}(m-m_0),\\
  \label{e:second_var}
  J''(m) &= \E^{\nu_0}[\Phi''_\mu(v+m)] + C_0^{-1}.
\end{align}

\subsection{Application of Robbins-Monro}

In \cite{Pinski:2015jy}, it was suggested that rather than try to find
a root of \eqref{e:el1}, the equation first be preconditioned by
multiplying by $C_0$,
\begin{equation}
  \label{e:el2}
  C_0\E^{\nu_0}[\Phi'_\mu(v+m)] + (m-m_0),
\end{equation}
and a root of this mapping is sought, instead.  Defining
\begin{subequations}
  \begin{align}
    f(m) &= C_0\E^{\nu_0}[\Phi'_\mu(v+m)] + (m-m_0),\\
    F(m,v) &= C_0\Phi'_\mu(v+m) + (m-m_0).
  \end{align}
\end{subequations}
The Robbins-Monro formulation is then
\begin{equation}
  m_{n+1} = m_n - a_{n+1} F(m_n, v_{n+1}) + P_{n+1},
\end{equation}
with $v_n \sim \nu_0$, i.i.d.

We thus have
\begin{thm}
  Assume:
  \begin{itemize}
  \item There exists $\nu=N(m, C_0)\sim \mu_0$ such that
    $\nu\ll\mu$.
  \item $\Phi_\mu'$ and $\Phi_\mu''$ exist for all $u \in \calH^1$.
  \item There exists $m_\star$, a local minimizer of $J$, such that
    $J'(m_\star)=0$.

  \item The mapping
    \begin{equation}
      \label{e:bounded1}
      m\mapsto \E^{\nu_0}\bracket{\norm{\sqrt{C_0}\Phi_\mu'(m+v)}^2}
    \end{equation}
    is bounded on bounded subsets of $\calH^1$.
  \item There exists a convex neighborhood $U_\star$ of $m_\star$ and
    a constant $\alpha>0$, such that for all $m\in U_\star$, for all
    $u \in \calH^1$,
    \begin{equation}
      \label{e:convexity1}
      \inner{J''(m)u}{u}\geq \alpha \norm{u}_{\calH^1}^2
    \end{equation}
  \end{itemize}
  Then, choosing $a_n$ according to Assumption \ref{a:sequence},
  \begin{itemize}
  \item If the subset $U_\star$ can be taken to be all of $\calH^1$,
    for the expanding truncation algorithm, $m_n \to m_\star$ a.s. in
    $\calH^1$.
  \item If the subset $U_\star$ is not all of $\calH^1$, then, taking
    $U_1$ to be a bounded (in $\calH^1$) convex subset of $U_\star$, with
    $m_\star \in U_1$, and $U_0$ any subset of $U_1$ such that there
    exist $R_0<R_1$ with
    \[
    U_0 \subset B_{R_0}(x_\star) \subset B_{R_1}(x_\star)\subset U_1,
    \]
    for the fixed truncation algorithm, $m_n\to m_\star$ a.s. in
    $\calH^1$.
  \end{itemize}
\end{thm}

\begin{proof}
We split the proof into 2 steps:
\begin{enumerate}[\upshape 1.]
    \item   By the assumptions of the theorem, we clearly satisfy Assumptions
  \ref{a:root} and \ref{a:sequence}.  To satisfy Assumption
  \ref{a:bound}, we observe that
  \begin{equation*}
    \E^{\nu_0}[\norm{F(m,v)}^2_{\calH^1}]\leq
    2\E^{\nu_0}\bracket{\norm{\sqrt{C_0}\Phi_\mu'(m+v)}^2} + 2\norm{m-m_0}_{\calH^1}^2.
  \end{equation*}
  This is bounded on bounded subsets of
  $\calH^1$.
  
  \item Per the convexity assumption, \eqref{e:convexity1}, implies
  Assumption \ref{a:convex}, since, by the mean value theorem in
  function spaces,
  \begin{equation*}
    \begin{split}
      \inner{m-m_\star}{f(m)}_{\calH^1} &=
      \inner{m-m_\star}{C_0\bracket{f(m_\star)
          +f'(\tilde m)(m-m_\star) }}_{\calH^1}\\
      &= \inner{m-m_\star}{J''(\tilde m)(m-m_\star)}\geq
      \alpha\norm{m-m_\star}_{\calH^1}^2
    \end{split}
  \end{equation*}
  where $\tilde m$ is some intermediate point between $m$ and
  $m_\star$.  This completes the proof.
\end{enumerate}

\end{proof}

While condition \eqref{e:convexity1} is sufficient to obtain
convexity, other conditions are possible.  For instance, suppose there
is a convex open set $U_\star$ containing $m_\star$ and constant
$\theta\in [0,1)$, such that for all $m \in U_\star$,
\begin{equation}
  \label{e:convexity2}
  \inf_{\substack{u\in \calH\\ u\neq 0}} \frac{\inner{\E^{\nu_0}[\Phi''_\mu(v+m)]u}{u}}{\norm{u}^2}\geq -\theta\lambda_1^{-1},
\end{equation}
where $\lambda_1$ is the principal eigenvalue of $C_0$.  Then this
would also imply Assumption \ref{a:convex}, since
\begin{equation*}
  \begin{split}
    \inner{m-m_\star}{f(m)}_{\calH^1} &=
    \inner{m-m_\star}{C_0\bracket{f(m_\star)
        +f'(\tilde m)(m-m_\star) }}_{\calH^1}\\
    &= \inner{m-m_\star}{J''(\tilde m)(m-m_\star)}\\
    &\geq \norm{m-m_\star}_{\calH^1}^2 +
    \inner{m-m_\star}{\E^{\nu_0}[\Phi''_\mu(v+\tilde m)] (m-m_\star)}\\
    &\geq \norm{m-m_\star}_{\calH^1}^2 -\theta \lambda_1^{-1} \norm{m-m_\star}^2\\
    &\geq (1-\theta)\norm{m-m_\star}_{\calH^1}^2.
  \end{split}
\end{equation*}
We mention \eqref{e:convexity2} as there may be cases, shown below,
for which the operator $\E^{\nu_0}[\Phi''_\mu(v+ m)]$ is obviously
nonnegative.

\section{Examples}
\label{s:examples}

To apply the Robbins-Monro algorithm to the relative entropy
minimization problem, the $\Phi_\mu$ functional of
interest must be examined.  In this section we present a few examples,
based on those presented in \cite{Pinski:2015jy}, and examine when the
assumptions hold.  The one outstanding assumption that we must make is
that, {\it a priori}, $\mu_0$ is an equivalent measure to $\mu$.

\subsection{Scalar Problem}
Taking $\mu_0 = N(0,1)$, the standard unit Gaussian, let $V:\R \to \R$
be a smooth function such that
\begin{equation}
  \frac{d\mu}{d\mu_0} = \frac{1}{Z_\mu} \exp\set{-{\eps^{-1}}V(x)}
\end{equation}
is a probability measure on $\R$.  In the above framework,
  \begin{align*}
    F(m,v) &= {\eps^{-1}}V'(m+v) +m ,\\
    f(m)  &= {\eps^{-1}}\E[V'(m+v)]+m  \\
    \Phi_\mu'(m) &= {\eps^{-1}}V'(m), \\
    \Phi_\mu''(m)& ={\eps^{-1}}V''(m)
  \end{align*}
and $v\sim N(0,1) = \nu_0 = \mu_0$.

\subsubsection{Globally Convex Case}
\begin{figure}
  {\includegraphics[width=6.25cm]{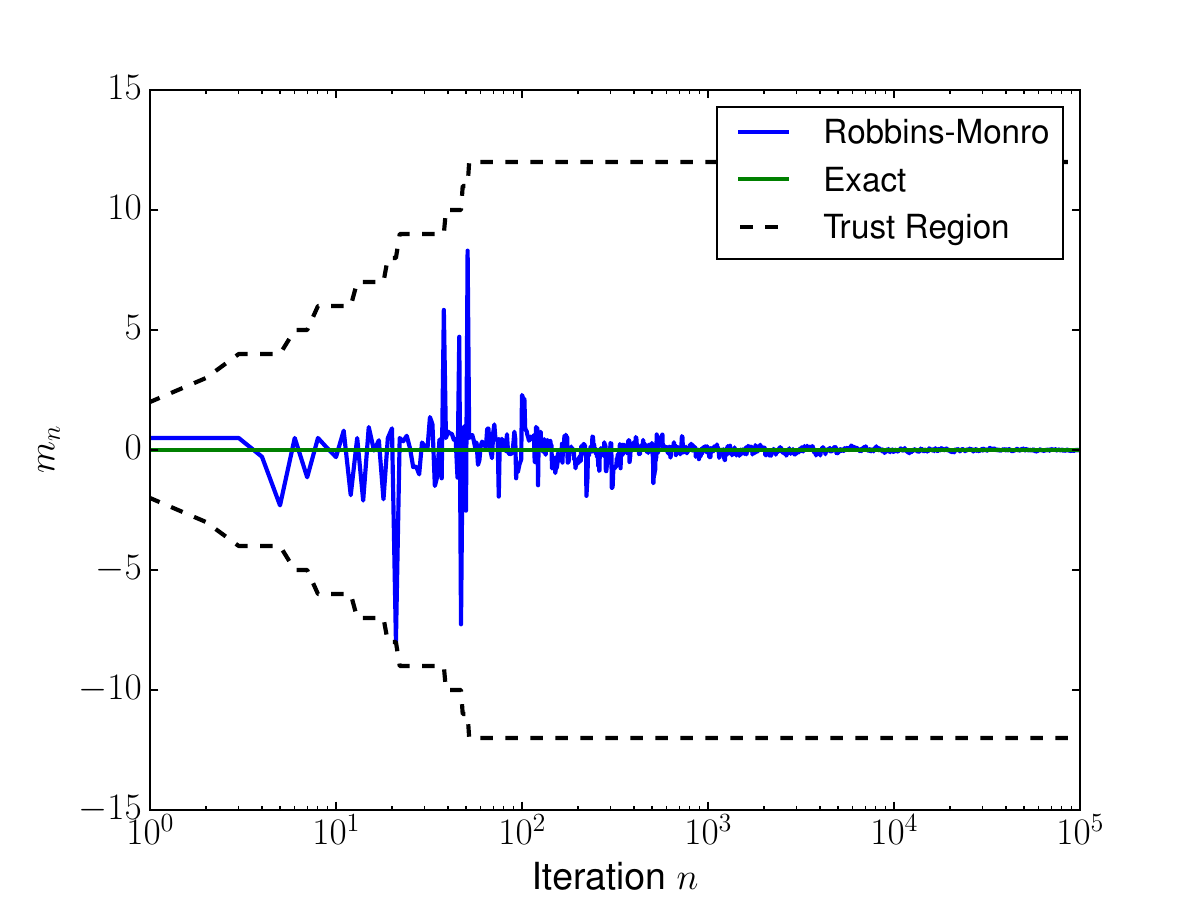}}
  \caption{Robbins-Monro applied to a globally convex scalar problem
    associated with \eqref{e:quartic} with $\eps=0.1$ and expanding
    trust regions $U_{n} = (-1-n,1+n)$. }
  \label{f:global_scalar}
\end{figure}

Consider the case that
\begin{equation}
  \label{e:quartic}
  V(x) = \tfrac{1}{2}x^2 + \tfrac{1}{4}x^4.
\end{equation}
In this case
  \begin{align*}
  F(m,v)& ={\eps^{-1}}\paren{m+v + (m+v)^3} + m,\\
  f(m) & = {\eps^{-1}}\paren{4m + m^3} +m, \\
  \E[\Phi''_\mu(m+v)] &= {\eps^{-1}}(4 + 3m^2),\\
  \E[\abs{\Phi'_\mu(m+v)}^2] &={\eps^{-1}} \paren{22 + 58 m^2 + 17
    m^4 + m^6}.
  \end{align*}
Since $\E[\Phi''_\mu(m+v)] \geq 4 {\eps^{-1}}$, all of our
assumptions are satisfied and the expanding truncation algorithm will
converge to the unique root at $m_\star =0$ a.s.  See Figure
\ref{f:global_scalar} for an example of the convergence at $\eps=0.1$,
$U_{n} = (-n -1,n+1)$, and always restarting at $0.5$.

We refer to this as a ``globally convex'' problem since $\calR$ is
globally convex about the minimizer.

\subsubsection{Locally Convex Case}

In contrast to the above problem, some minimizers are only ``locally'' convex.  Consider the case the double well potential
\begin{equation}
  \label{e:dblewell}
  V(x) = \tfrac{1}{4}(4-x^2)^2
\end{equation}
Now, the expressions for RM are
  \begin{align*}
  F(m,v) &= {\eps^{-1}}\paren{(m+v)^3-4(x+v))} + m,\\ 
  f(m)  &= {\eps^{-1}}\paren{m^3-x} +m, \\
  \E[\Phi''_\mu(m+v)] &={\eps^{-1}}\paren{3m^2-1},\\
  \E[\abs{\Phi'_\mu(m+v)}^2] &={\eps^{-1}} (1 + m^2) (7 + 6 m^2 +m^4).
  \end{align*}
In this case, $f(m)$ vanishes at $0$ and $\pm \sqrt{1-\eps}$, and
$J''$ changes sign from positive to negative when $m$ enters
$({-\sqrt{(1-\eps)/3},\sqrt{({1-\eps})/{3}}})$.  We must therefore
restrict to a fixed trust region if we want to ensure convergence to
either of $\pm\sqrt{1-\eps}$.

\begin{figure}
  \subfigure[]{\includegraphics[width=6.25cm]{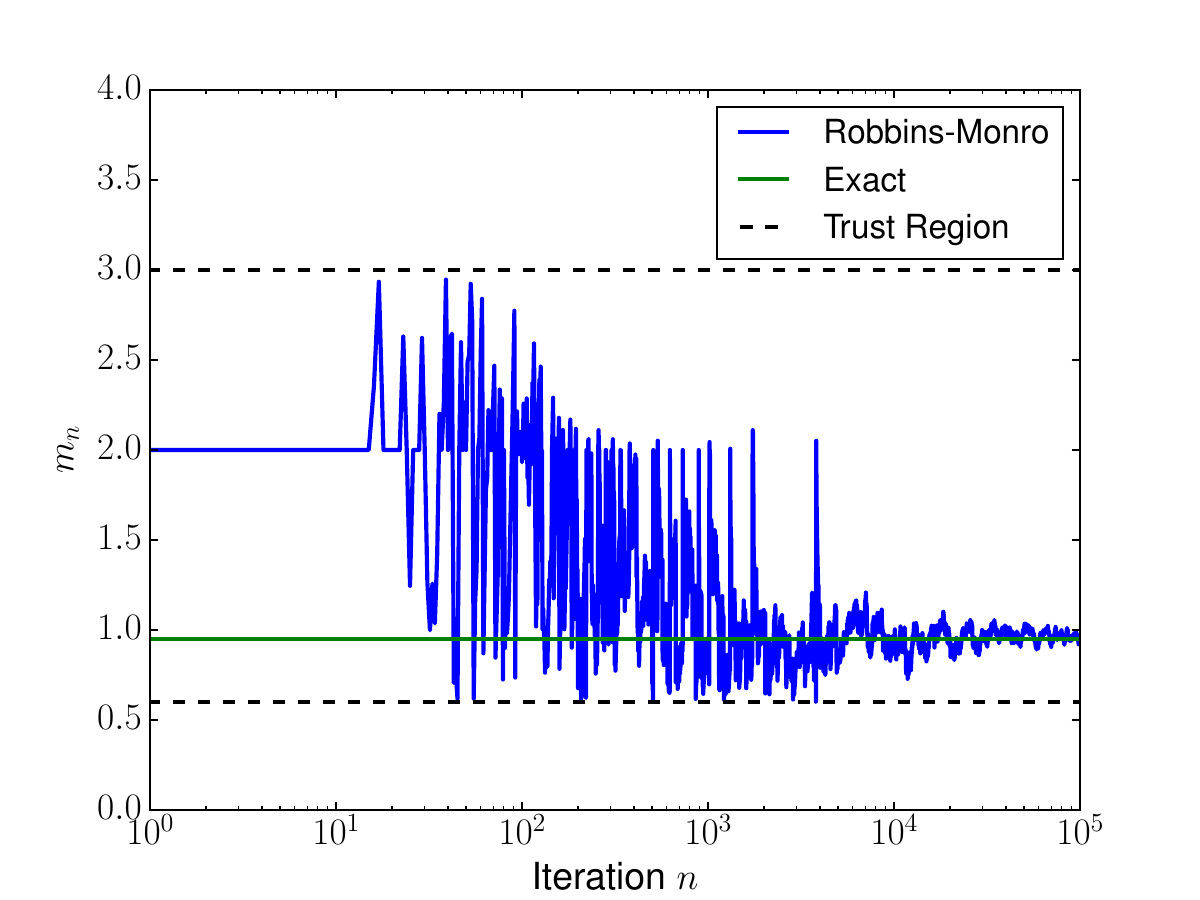}}
  \subfigure[]{\includegraphics[width=6.25cm]{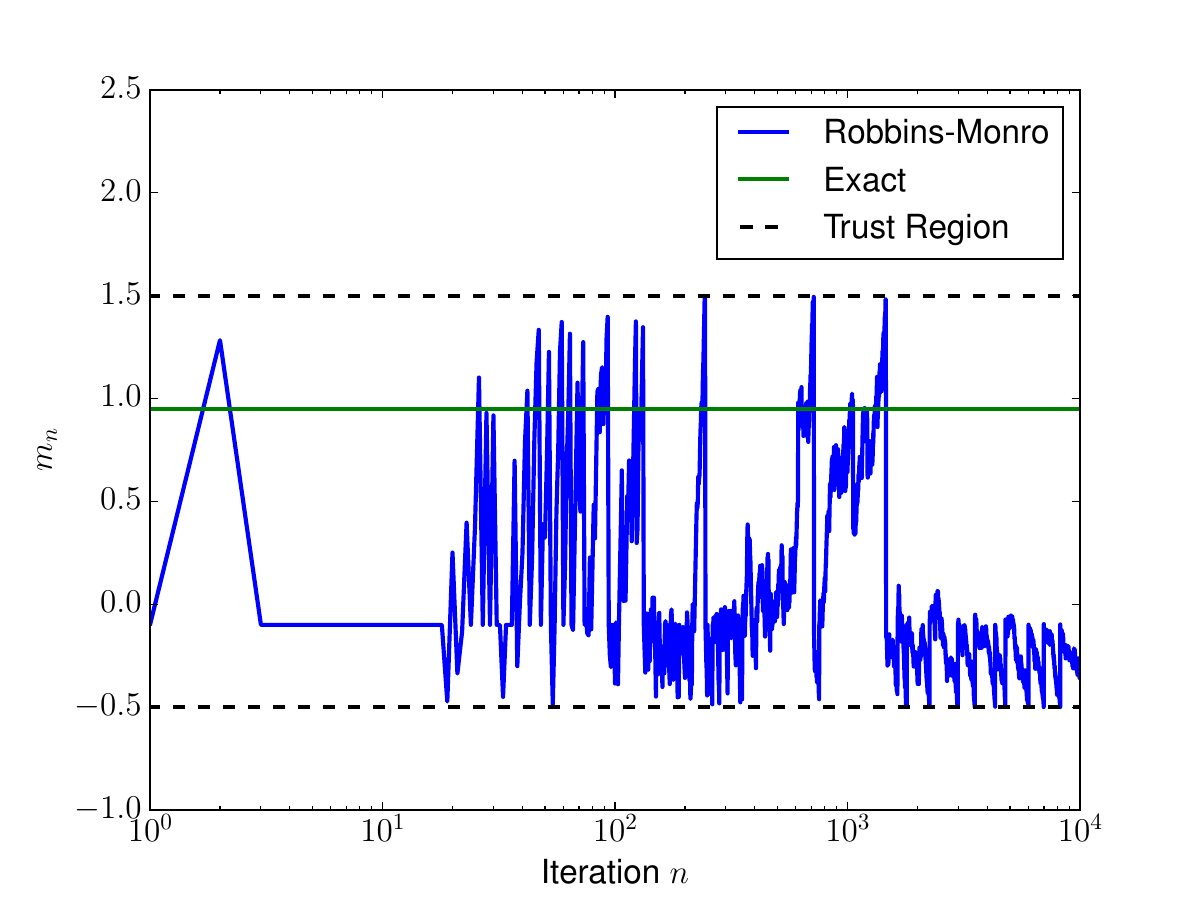}}
  \caption{Robbins-Monro applied to the nonconvex scalar problem
    associated with \eqref{e:dblewell}.  Figure (a) shows the result
    with a well chosen trust region, while (b) shows the outcome of a
    poorly chosen trust region.}
  \label{f:scalar}
\end{figure}

We ran the problem at $\eps = 0.1$ in two cases.  In the first case,
$U_1 = (0.6, 3.0)$ and the process always restarts at $2$.  This
guarantees convergence since the second variation will be strictly
positive.  In the second case, $U_1 = (-0.5, 1.5)$, and the process
always restarts at $-0.1$.  Now, the second variation
can change sign.  The results of these two experiments appear in
Figure \ref{f:scalar}.  For some random number sequences the algorithm still converged to $\sqrt{1-\eps}$, even with the
poor choice of trust region.

\subsection{Path Space Problem}
Take $\mu_0=N(m_0(t), C_0)$, with
\begin{equation}
  C_0 = \paren{-\frac{d^2}{dt^2}}^{-1},
\end{equation}
equipped with Dirichlet boundary conditions on
$\calH = L^2(0,1)$.\footnote{This is the covariance of the
  standard unit Brownian bridge, $Y_t = B_t - t B_1$.}   In this case
the Cameron-Martin space $\calH^1 = H^1_0(0,1)$, the standard Sobolev
space equipped with the Dirichlet norm.  Let us assume
$m_0 \in H^1(0,1)$, taking values in $\R^d$.

Consider the path space distribution on $L^2(0,1)$,
induced by
\begin{equation}
  \label{e:path_potential}
  \frac{d\mu}{d\mu_0} = - \frac{1}{Z_\mu}\exp\set{-\Phi_\mu(v)}, \quad
  \Phi_\mu(u) = {\eps^{-1}}\int_0^1 V(v(t))dt,
\end{equation}
where $V:\R^d\to\R$ is a smooth function.  We assume that $V$ is such
that this probability distribution exists and that $\mu \sim \mu_0$,
our reference measure.

We thus seek an $\R^d$ valued function $m(t) \in H^1(0,1)$ for our
Gaussian approximation of $\mu$, satisfying the boundary conditions
\begin{equation}
  \label{e:bc}
  m(0) = m_-,\quad m(1) = m_+.
\end{equation}
For simplicity, take $m_0 = (1-t)m_- + t m_+$, the linear interpolant
between $m_\pm$.  As above, we work in the shifted coordinated
$x(t) = m(t) - m_0(t)\in H^1_0(0,1)$.

Given a path $v(t)\in H^1_0$, by the Sobolev embedding,
$v$ is continuous with its $L^\infty$ norm controlled by its $H^1$
norm.  Also recall that for $\xi \sim N(0, C_0)$, in the case of
$\xi(t) \in \R$,
\begin{equation}
  \label{e:bb_moments}
  \E\bracket{\xi(t)^p} = \begin{cases}
    0, & \text{$p$ odd},\\
    (p-1)!!\bracket{t(1-t)}^{\frac p 2}, & \text{$p$ even}.
  \end{cases}
\end{equation}
Letting $\lambda_1 =1/\pi^2$ be the ground state eigenvalue of $C_0$,
\begin{equation*}
  \begin{split}
    \E[\|\sqrt{C_0}\Phi'_\mu(v +m_0+\xi)\|^2]&\leq
    {\lambda_1}\E[\|\Phi'_\mu(v +m_0+\xi)\|^2]\\
    &\quad= {\lambda_1}{\eps^{-2}}\int_0^1
    \E[{\abs{V'(v(t)+m_0(t)+\xi(t))}^2}]dt.
  \end{split}
\end{equation*}
The terms involving $v+m_0$ in the integrand can be controlled by the
$L^\infty$ norm, which in turn is controlled by the $H^1$ norm, while
the terms involving $\xi$ can be integrated according to
\eqref{e:bb_moments}.  As a mapping applied to $x$, this expression is
bounded on bounded subsets of $H^1$.

Minimizers will satisfy the ODE
\begin{equation}
  \label{e:ode}
  {\eps}^{-1}\E\bracket{V'(x+m_0 +\xi)} -x'' = 0,\quad x(0) =x(1)=0.
\end{equation}

\subsection{Globally Convex Example}

\begin{figure}
  \subfigure[]{\includegraphics[width=6.25cm]{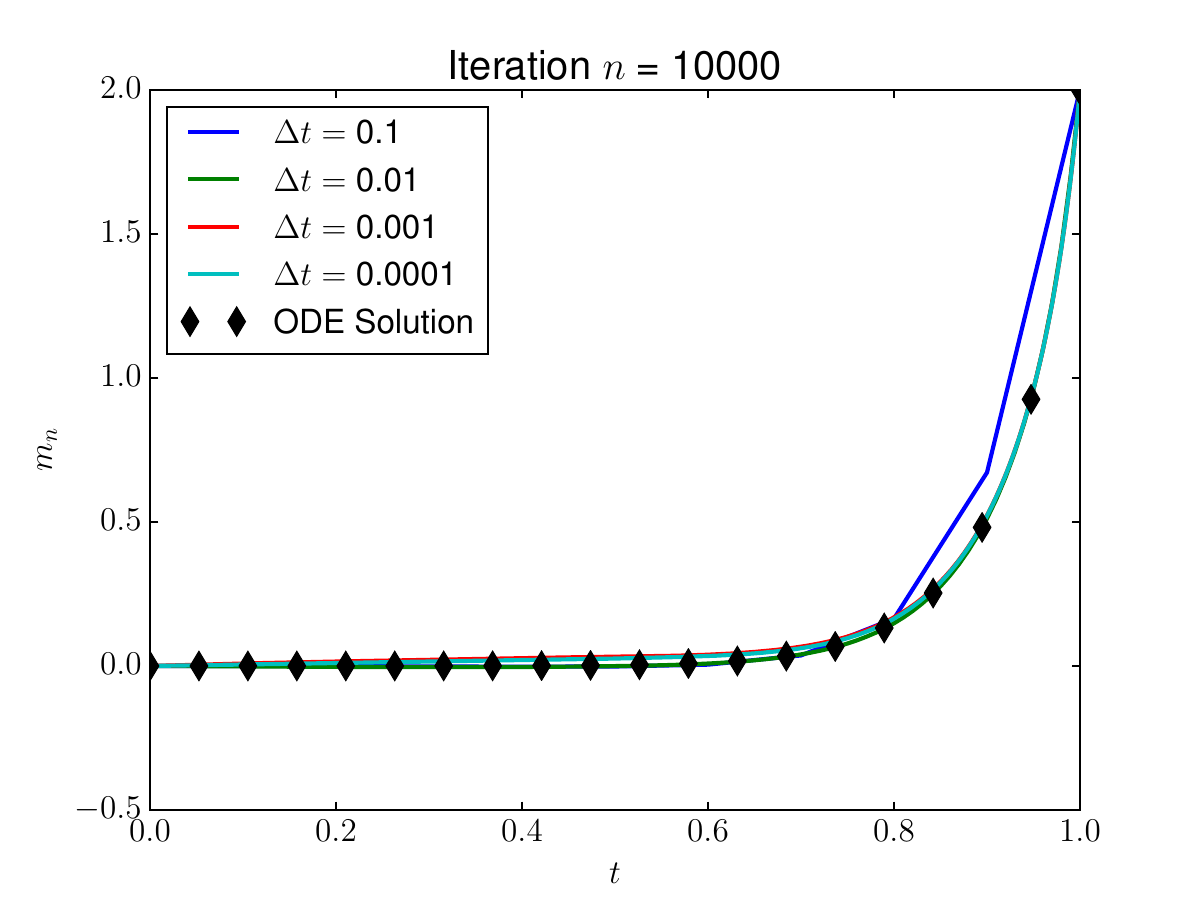}}
  \subfigure[]{\includegraphics[width=6.25cm]{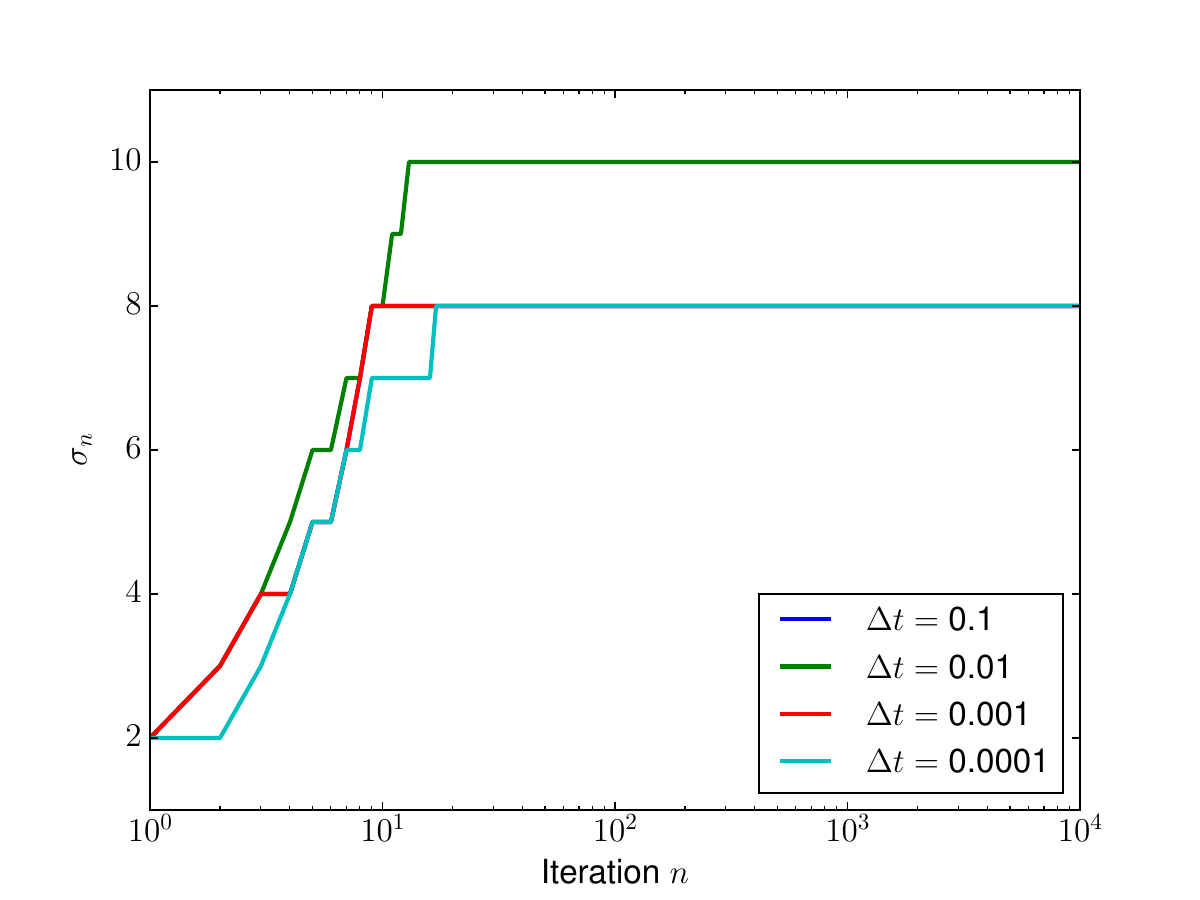}}
  \caption{The mean paths computed for \eqref{e:quartic_path} at
    different resolutions, along with the truncation sequence.}
  \label{f:quartic_path}
\end{figure}

With regard to convexity about a minimizer, $m_\star$, if, for
instance, $V''$ were  pointwise positive definite, then the
problem would satisfy \eqref{e:convexity2}, ensuring convergence.
Consider the quartic potential $V$ given by \eqref{e:quartic}.  In
this case,
\begin{equation}
  \label{e:quartic_path}
  \Phi(v) = {\eps}^{-1}\int_0^1 \frac{1}{2}v(t)^2 +\frac{1}{4}v(t)^4 dt,
\end{equation}
and
\begin{align*}
  \Phi'(v+m_0+ \xi) & = {\eps}^{-1}\bracket{(v+m_0 +\xi) +3(v+m_0 +\xi)^3 }\\
  \Phi''(v+m_0 + \xi) & = {\eps}^{-1}\bracket{1 +3 (v+m_0+\xi)^2},\\
  \E[\Phi'(v+m_0 + \xi)]& = {\eps}^{-1}\bracket{v+m_0 +(v+m_0)^3+
                          3  t(1-t) (v+m_0)}\\
  \E[\Phi''(v+m_0 + \xi)]& = {\eps}^{-1}\bracket{1 + 3 (v+m_0)^2 +
                           3  t(1-t) }
\end{align*}
Since $ \Phi''(v+m_0+\xi)\geq \eps^{-1}$, we are guaranteed convergence
using expanding trust regions.  Taking $\eps=0.01$, $m_- = 0$ and
$m_+ = 2$, this is illustrated in Figure \ref{f:quartic_path}, where
we have also solved \eqref{e:ode} by ODE methods for comparison.  As
trust regions, we take
\begin{equation}
  U_n = \set{m \in H^1_0(0,1)\mid \norm{x}_{H^1}\leq 10+n},
\end{equation}
and we always restart at the zero solution
Figure \ref{f:quartic_path} also shows robustness to discretization;
the number of truncations is relatively insensitive to $\Delta t$.

\subsection{Locally Convex Example}
For many problems of interest, we do not have global convexity.
Consider the double well potential \eqref{e:dblewell}, but in the case
of paths,
\begin{equation}
  \label{e:dblewell_path}
  \Phi(u) = {\eps^{-1}}\int_0^1\frac{1}{4} (4-v(t)^2)^2dt.
\end{equation}
Then,
\begin{align*}
  \Phi'(v + m_0 + \xi)& = {\eps}^{-1}\bracket{(v + m_0 + \xi)^3 - 4 (v + m_0 + \xi)}\\
  \Phi''(v + m_0 + \xi) & = {\eps}^{-1}\bracket{3 (v + m_0 + \xi)^2 - 4},\\
  \E[\Phi'(v+m_0 + \xi)]& = {\eps}^{-1}\bracket{(v+m_0)^3 +
                          3  t(1-t) (v+m_0)-4(v+m_0)}\\
  \E[\Phi''(v+m_0 + \xi)]& = {\eps}^{-1}\bracket{3(v+m_0)^2 +
                           3  t(1-t) -4}
\end{align*}

\begin{figure}
  {\includegraphics[width=6.25cm]{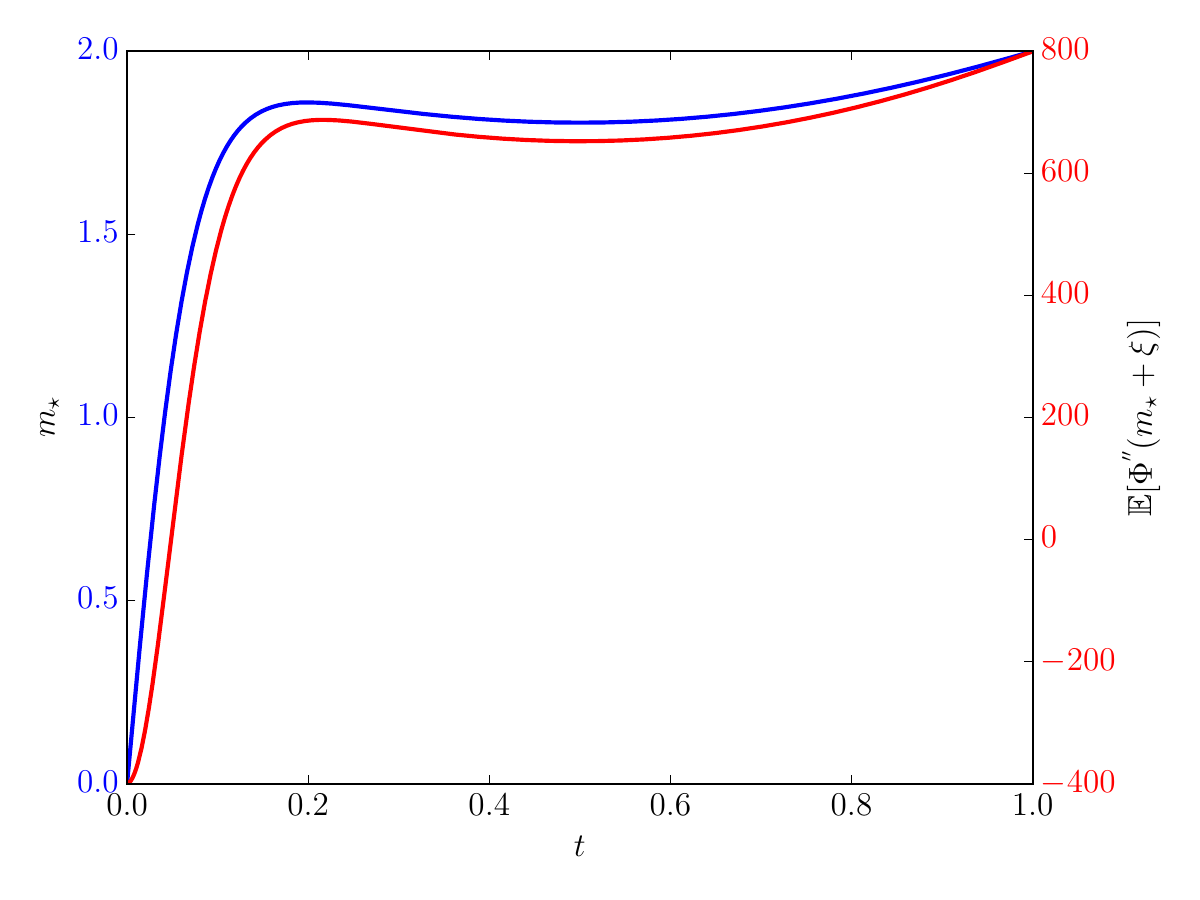}}
  \caption{The numerically computed solution to \eqref{e:ode} in the
    case of the double well, \eqref{e:dblewell_path}, $m_\star$, and
    the associated $\E^{\nu_0}[\Phi''(m_\star + \xi)]$.}
  \label{f:path1}
\end{figure}

Here, we take $m_-= 0$, $m_+ = 2$, and $\eps = 0.01$.  We
have plotted the numerically solved ODE in Figure \ref{f:path1}.  Also
plotted is $\E[\Phi''(v_\star +m_0+ \xi)]$.  Note that $\E[\Phi''(v_\star +m_0+ \xi)]$ is not sign definite,
becoming as small as $-400$.  Since $C_0$ has
$\lambda_1 = 1/\pi^2 \approx 0.101$, \eqref{e:convexity2} cannot
apply.

\begin{figure}
  {\includegraphics[width=6.25cm]{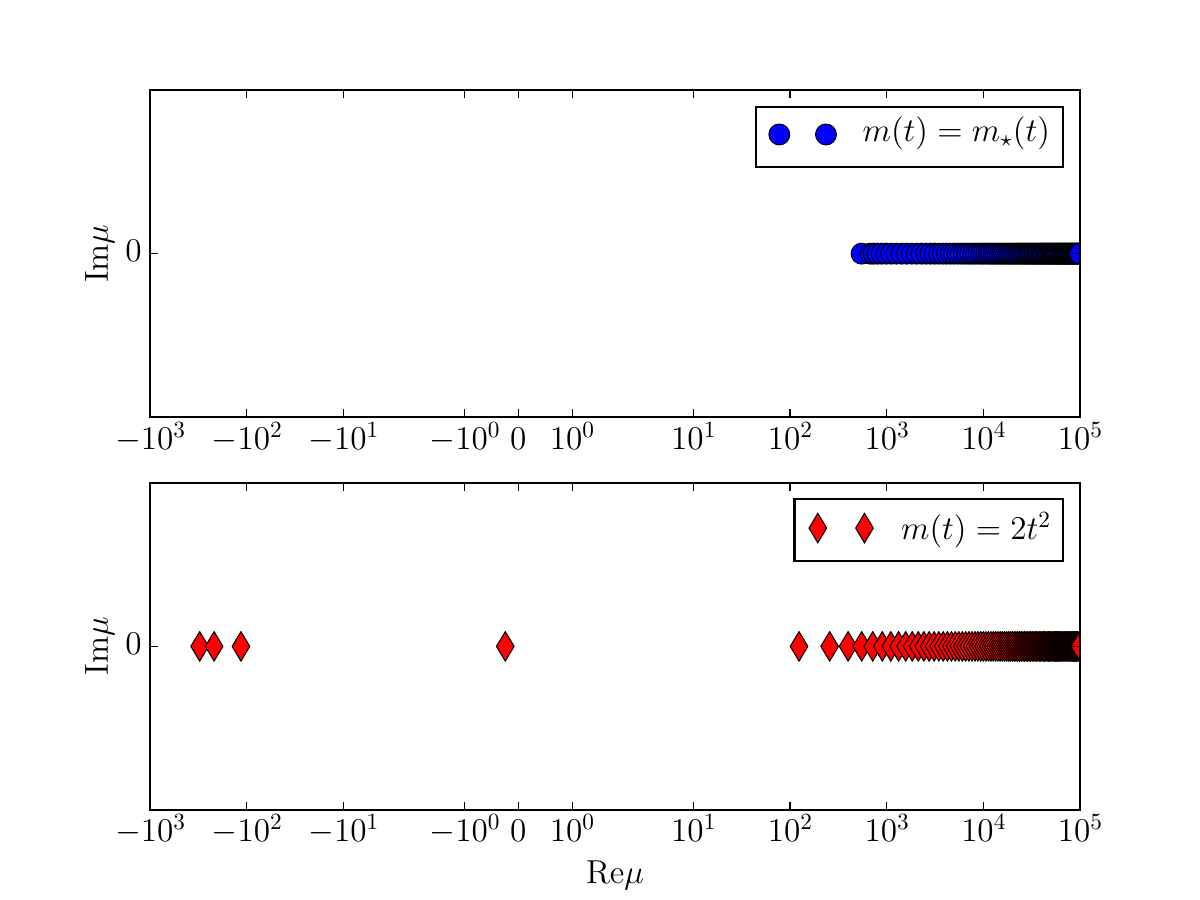}}
  \caption{The numerically computed spectrum for \eqref{e:soperator},
    associated with the $m_\star$ shown in Figure \ref{f:path1}.  Also
    shown is the numerically computed spectrum for the path
    $m(t) = 2t^2$, which introduces negative eigenvalues.}
  \label{f:spec}
\end{figure}

Discretizing the Schr\"odinger operator
\begin{equation}
  \label{e:soperator}
  J''(v_\star) = -\frac{d^2}{dt^2} + {\eps}^{-1}\paren{3(v_\star(t)+m_0(t))^2 +
    3  t(1-t) -4},
\end{equation}
we numerically compute the eigenvalues.  Plotted in Figure
\ref{f:spec}, we see that the minimal eigenvalue of $J''(m_\star)$ is
approximately $\mu_1\approx 550$.  Therefore,
\begin{equation}
  \inner{J''(x_\star)u}{u}\geq \mu_1 \norm{u}^2_{L^2}\Rightarrow \inner{J''(x)u}{u}\geq \alpha\norm{u}_{H^1}^2,
\end{equation}
for all $v$ in some neighborhood of $v_\star$.  For an appropriately
selected fixed trust region, the algorithm will converge.

However, we can show that the convexity condition is not global.
Consider the path $m(t) = 2t^2$, which satisfies the boundary
conditions .  As shown in Figure \ref{f:spec}, this path induces
negative eigenvalues.

\begin{figure}
  \subfigure[]{\includegraphics[width=6.25cm]{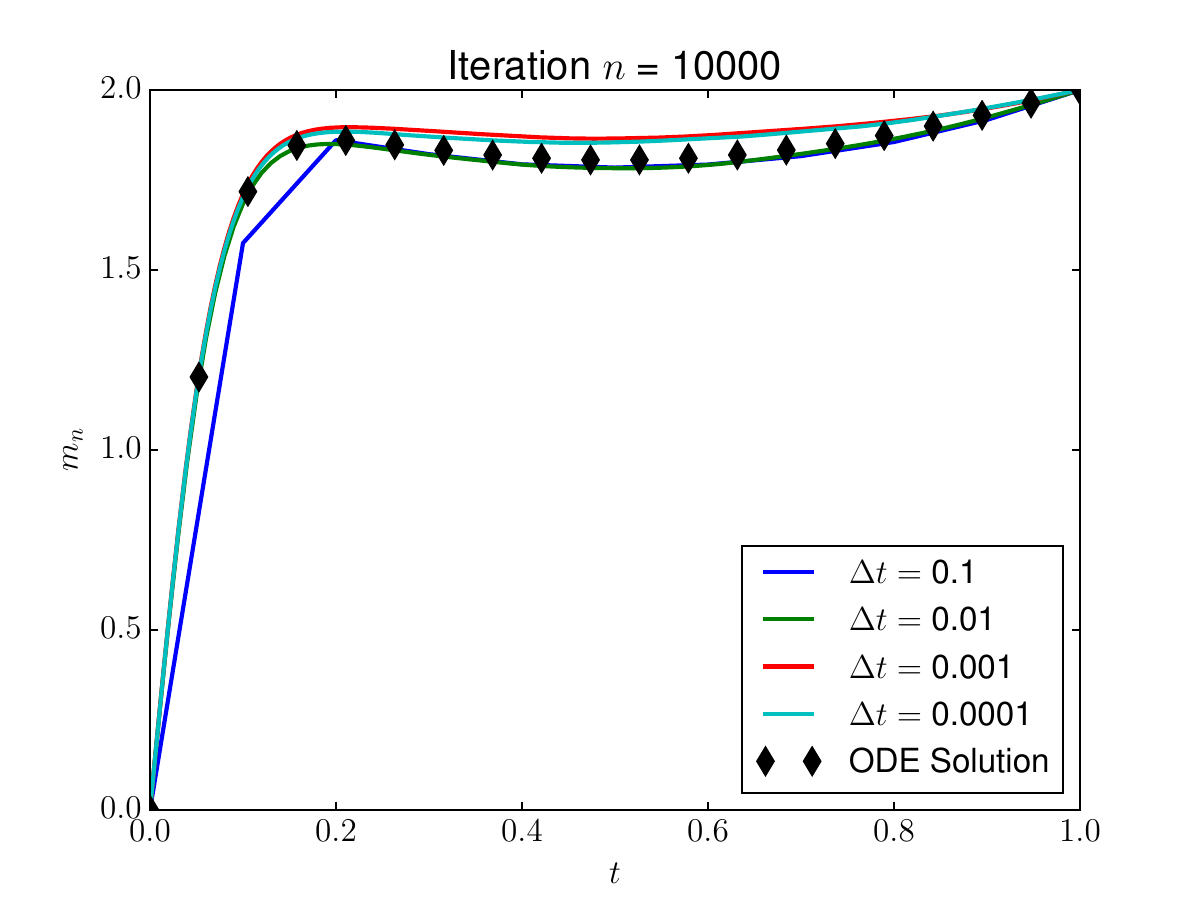}}
  \subfigure[]{\includegraphics[width=6.25cm]{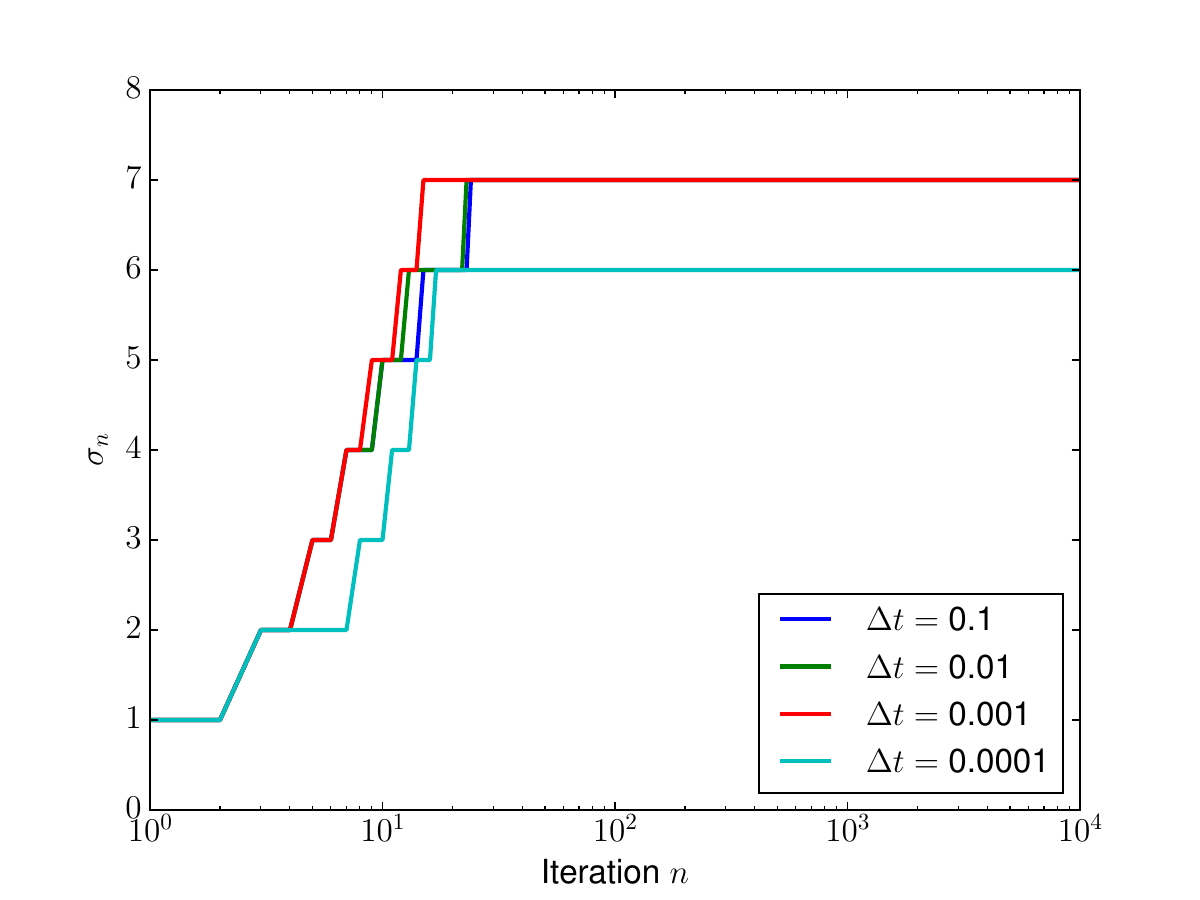}}
  \caption{The mean paths computed for \eqref{e:dblewell_path} at
    different resolutions, along with the truncation sequence.}
  \label{f:doublewell_path}
\end{figure}

Despite this, we are still observe convergence.  Using the fixed trust
region
\begin{equation}
  \label{e:path_trust2}
  U_1 = \set{x\in H^1_0(0,1)\mid \norm{x}_{H^1}\leq 100},
\end{equation}
we obtain the results in Figure \ref{f:doublewell_path}.  Again,
the convergence is robust to discretization.

\section{Discussion}
\label{s:disc}
We have shown that the Robbins-Monro algorithm, with both fixed and
expanding trust regions, can be applied to Hilbert space valued
problems, adapting the finite dimensional proof of
\cite{Lelong:2008ck}.  We have also constructed sufficient conditions
for which the relative entropy minimization problem fits within this
framework.

One problem we did not address here was how to identify fixed trust
regions.  Indeed, that requires a tremendous amount of {\it a priori}
information that is almost certainly not available.  We interpret that
result as a local convergence result that gives a theoretical basis
for applying the algorithm.  In practice, since the root is likely
unknown, one might run some numerical experiments to identify a
reasonable trust region, or just use expanding trust regions.  The
practitioner will find that the algorithm converges to a solution,
though perhaps not the one originally envisioned.  A more
sophisticated analysis may address the convergence to a set of roots,
while being agnostic as to which zero is found.

Another problem we did not address was how to optimize not just the
mean, but also the covariance in the Gaussian.  As discussed in
\cite{Pinski:2015jy}, it is necessary to parameterize the covariance
in some way, which will be application specific.  Thus, while the form
of the first variation of relative entropy with respect to the mean,
\eqref{e:el1}, is quite generic, the corresponding expression for the
covariance will be specific to the covariance parameterization.
Additional constraints are also necessary to guarantee that the
parameters always induce a covariance operator.  We leave such
specialization as future work.

\section*{Acknowledgments}

This work was supported by US Department of Energy Award DE-SC0012733.  This
work was completed under US National Science Foundation Grant DMS-1818716.  The
authors would like to thank J. Lelong for helpful comments, along with anonymous reviewers whose reports significantly impacted our work.

\bibliographystyle{plain}

\bibliography{rm_refs}

\end{document}